\documentclass[11pt]{article}

\usepackage[compress]{natbib}
\usepackage{amsmath,amssymb,amsthm}
\usepackage{color}
\usepackage[margin=1in]{geometry}
\usepackage{setspace}
\usepackage{amsthm, amsmath, amssymb, amsfonts}
\usepackage{subcaption}
\usepackage{bm}
\usepackage[dvipsnames]{xcolor}
\usepackage{hyperref}
\hypersetup{pdftex, colorlinks=true, citecolor=MidnightBlue, linkcolor=red!80!black, urlcolor=MidnightBlue}
\usepackage{caption}
\usepackage{url}            
\usepackage{booktabs}       
\usepackage{nicefrac}       
\usepackage{graphicx}
\usepackage{caption}
\usepackage{thm-restate}
\usepackage[ruled,linesnumbered]{algorithm2e}
\usepackage{algpseudocode}
\usepackage{bbm}
\usepackage[export]{adjustbox}

\usepackage{microtype}
\usepackage{mathpazo} 
\usepackage[scaled]{helvet} 
\usepackage{courier} 
\normalfont
\usepackage[T1]{fontenc}

\usepackage{threeparttable}
\allowdisplaybreaks[4] 
\usepackage[shortlabels]{enumitem} 

\usepackage{tikz}
\usetikzlibrary{fit}
\tikzset{%
	highlight/.style={rectangle,blend mode = multiply,draw=blue!90!black,thick,rounded corners = 0.3 mm,inner sep=0.5pt}
}

\usepackage{tcolorbox}
\def\O#1{\text{\ding{\the\numexpr#1+171}}}

\newcommand{\ceil}[1]{\left\lceil #1 \right\rceil}

\newcommand{\bz}{\mathbf 0}

\newcommand{\R}{\mathbb R}

\DeclareMathOperator{\prox}{prox}
\DeclareMathOperator{\proj}{proj}

\DeclareMathOperator{\dist}{dist}
\DeclareMathOperator*{\argmin}{argmin}

\def\B{{\mathbb{B}}}

\def\N{{\mathbb{N}}}
\def\O{{\mathbb{O}}}

\def\R{{\mathbb{R}}}

\def\bA{{\mathbf{A}}}

\def\cB{{\cal B}}
\def\cC{{\cal C}}
\def\cD{{\cal D}}

\def\cI{{\cal I}}
\def\cJ{{\cal J}}

\def\cM{{\cal M}}

\def\cO{{\cal O}}

\def\cR{{\cal \cR}}

\def\cX{{\cal X}}

\def\b{\bm{b}}

\def\e{\bm{e}}

\def\p{\bm{p}}

\def\u{\bm{u}}
\def\v{\bm{v}}

\def\x{\bm{x}}
\def\y{\bm{y}}
\def\z{\bm{z}}
\def\bz{{\mathbf 0}}

\def\1{{\mathbf 1}}
\def\intm{{\rm int}}
\def\bmu{\bm{\mu}}
\def\blam{\bm{\lambda}}

\DeclareMathOperator{\cl}{cl}
\DeclareMathOperator{\dom}{dom}

\DeclareMathOperator{\rank}{rank}

\DeclareMathOperator{\inter}{int}

\DeclareMathOperator{\bd}{bd}

\usepackage{thmtools}

\declaretheoremstyle[parent=section]{definitionwithend}

\declaretheorem[style=definitionwithend]{theorem}
\declaretheorem[style=definitionwithend]{proposition}
\declaretheorem[style=definitionwithend]{definition}
\declaretheorem[style=definitionwithend]{assumption}
\declaretheorem[style=definitionwithend]{example}
\declaretheorem[style=definitionwithend]{remark}

\declaretheorem[style=definitionwithend]{lemma}
\declaretheorem[style=definitionwithend]{fact} 
\declaretheorem[style=definitionwithend]{observation}

\title{Spurious Stationarity and  Hardness Results for 
	Bregman Proximal-Type Algorithms \footnote{{Authors are listed in alphabetical order.}}}

\date{\today}
\author{%
	He Chen\thanks{Department of Systems Engineering and Engineering Management, The Chinese University of Hong Kong, Shatin, NT, Hong Kong. \texttt{hchen@se.cuhk.edu.hk}} \and
	Jiajin Li\thanks{Sauder School of Business, University of British Columbia, Vancouver, BC, Canada. \texttt{jiajin.li@sauder.ubc.ca}} \and
	Anthony Man-Cho So\thanks{Department of Systems Engineering and Engineering Management, The Chinese University of Hong Kong, Shatin, NT, Hong Kong. \texttt{manchoso@se.cuhk.edu.hk}}
}

\begin{document}
	\maketitle
	
	\begin{abstract}
Bregman proximal-type algorithms (BPs), such as mirror descent, are widely used for exploiting problem structures through non-Euclidean geometries. 
In this paper, we show that BPs can get trapped near a class of non-stationary points, which we term \emph{spurious stationary points}. 
Such stagnation can persist for any prescribed finite number of iterations when the gradient of the Bregman kernel diverges near the boundary of its domain, even in convex linear programs.
The underlying mechanism is a boundary-induced mismatch between the geometry of Bregman updates and the standard first-order stationarity condition. 
Near the boundary, a vanishing Bregman residual, measured by the Bregman
divergence between an iterate and its update, may certify only stationarity on a
lower-dimensional affine manifold induced by the kernel, rather than
stationarity of the original problem. Consequently, small Bregman residuals may
suggest convergence even when the iterates remain near non-stationary points.
We further show that this phenomenon is not merely an artifact of poor initialization: In a nonconvex example, Bregman proximal gradient (BPG) methods can be drawn toward a spurious stationary point from a well-behaved interior initialization, while Euclidean projected gradient dynamics moves toward the global minimizer. 
Finally, we prove that spurious stationary points are not isolated artifacts; for polyhedral feasible sets and kernels whose closed domains coincide with the
nonnegative orthant, every non-stationary vertex is spurious.
Together, our findings reveal a structural blind spot in Bregman residual-based convergence analysis and motivate alternative certificates or algorithmic safeguards.
	\end{abstract} 
	
	\section{Introduction}
	In this paper, we consider structured nonsmooth (non)convex optimization problems of the form 
	\begin{equation}\label{eq:obj}
		\min\limits_{\x \in\R^n} F(\x):= f(\x)+ g(\x),
		\tag{P}
	\end{equation}
	where $\dom(g)=\cX$ is a nonempty closed convex set, $f:\R^n \rightarrow  \R$ is a continuously differentiable function, and $ g:\R^n \rightarrow \overline \R$ is a convex and locally Lipschitz continuous function on $\cX$ relative to its domain. 

To solve \eqref{eq:obj}, Bregman proximal-type algorithms (BPs) 
are widely used for leveraging the geometry of $\cX$ while avoiding costly Euclidean projections or proximal operations; see, e.g., \citep{beck2003mirror, arora2012multiplicative, zhang2021proximal}.

In this work, we study BPs through the unified update	
	\begin{equation}\label{eq:unified}
\x^{k+1}=\mathop{\argmin}\limits_{\y \in \R^n}\ \left\{\gamma\left(\y;\x^k\right)+ g(\y)+\frac{1}{t_k}D_h(\y,\x^k)\right\}, 
	\end{equation}  
	where $\gamma(\ \cdot\ ;\x)$ is the surrogate model for $f$ at point $\x$, $t_k > 0$ is the step size, and  $D_h$ denotes the Bregman divergence induced by a kernel function $h$. Many classical algorithms fit within this framework. For example, setting $\gamma(\y;\x^k) =  f(\x^k)+\nabla f(\x^k)^T(\y-\x^k)$ recovers the \emph{Bregman proximal gradient method}  (BPG) \citep{censor1992proximal,bauschke2017descent,bauschke2019linear,zhu2021level}.  Choosing $\gamma = f$ gives the \emph{Bregman proximal point method} \citep{chen1993convergence,kiwiel1997proximal}. Moreover, using a second-order surrogate,  $\gamma(\y;\x^k) := f(\x^k)+\nabla f(\x^k)^T(\y-\x^k)+\frac12(\y-\x^k)^T\nabla^2f(\x^k)(\y-\x^k)$, leads to a second-order variant recently studied by \citet{doikov2023gradient}. 
    
The central finding of this paper is that BPs can become trapped near certain non-stationary points, which we term \emph{spurious stationary points}. When the gradient of the Bregman kernel $h$ diverges near the boundary of its domain, such trapping can persist for any prescribed finite number of iterations. This behavior stands in sharp contrast to Euclidean gradient methods, as well as to BPs equipped with either full-domain kernels, i.e., $\dom(h)=\R^n$, or kernels with Lipschitz continuous gradients; see, e.g., \citep{zhang2025stochastic,zhang2018convergence}. These regular settings preclude the boundary-induced geometric degeneration that gives rise to spurious stationarity in our analysis.

To illustrate this boundary-induced failure mode, we first present a simple
linear programming (LP) problem for which BPG with the entropy kernel can
remain near a non-optimal boundary point for an arbitrarily prescribed finite
number of iterations.

    \begin{quote}
\textbf{Illustrative Example -- Finite-Time Stagnation of BPG.}
Consider the following simple LP problem:
\begin{equation}
		\label{example:lp}
		\begin{array}{rl}
			\min\limits_{x_1,x_2}& -x_1 \\
			{\rm s.t.}  & x_1+x_2=1,x_1,x_2\geq0,
		\end{array} 
	\end{equation}
    which  admits the unique solution at $(1,0)$. If we choose the Boltzmann–Shannon entropy kernel, i.e.,   $h(\x) =\sum_{i=1}^2 x_i\log x_i$, BPG admits the closed-form iteration:  
\begin{align*}
		\x^{k+1} 
		& = \left(\frac{x_1^k}{x_1^k+e^{-t}x_2^k}, \frac{e^{-t}x_2^k}{x_1^k+e^{-t}x_2^k}\right), \quad \forall~ k \in \N_+.
	\end{align*} 
    Although $(1,0)$ is the unique solution, the non-optimal boundary point
$(0,1)$ behaves as a spurious stationary point under the entropy geometry.
Indeed, for any prescribed iteration budget $K\in\N$ and any tolerance
$\epsilon>0$, one can choose an interior feasible initialization sufficiently
close to $(0,1)$ such that every iterate $\x^k$, $0\leq k\leq K$, remains
within an $\epsilon$-neighborhood of $(0,1)$. Thus, even in this convex LP,
BPG can exhibit arbitrarily long finite-time stagnation near a non-stationary
boundary point. Figure~\ref{fig:enter-label} illustrates this finite-time stagnation phenomenon.
	\begin{figure}
		\centering
		\includegraphics[width=0.6\textwidth]{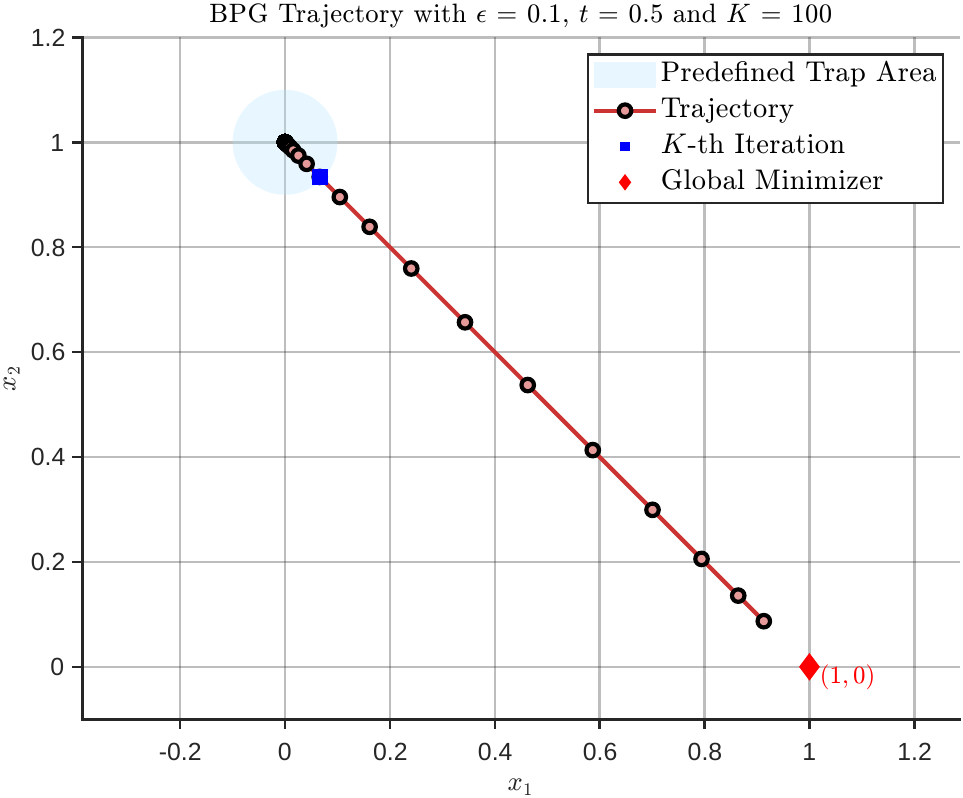}
		\caption{
        The trajectory of BPG with the Boltzmann--Shannon entropy kernel on the LP instance~\eqref{example:lp}.
        For a suitably chosen initialization, all iterates up to the prescribed horizon $K$ remain near the non-optimal boundary point $(0,1)$, while the unique global minimizer is $(1,0)$. 
        The shaded region indicates the predefined $\epsilon$-neighborhood.
            }
		\label{fig:enter-label}
	\end{figure}
\end{quote}
The example highlights a purely geometric source of finite-time stagnation:
Even in this elementary LP, the entropy-induced BPG update can become
arbitrarily small near the non-stationary boundary point $(0,1)$.
This raises the following question:
\begin{center}
\bf{
How can one rigorously define such spurious stationary points \\ and characterize
the finite-time stagnation of BPs near them?
}
\end{center}
In what follows, we formalize this phenomenon and present our main theoretical results. 


Our first result identifies the precise variational meaning of a vanishing Bregman residual near the boundary of the kernel domain.
For a fixed stepsize $t>0$, let $T_\gamma^t(\x)$ denote the one-step Bregman update in~\eqref{eq:unified} with $t_k=t$, and define $ R_\gamma^t(\x):=D_h(T_\gamma^t(\x),\x),$
a standard Bregman residual used as a stationarity measure in the analysis of BPs; see, e.g., \citep{bauschke2017descent, zhang2018convergence, bedi2022hidden, latafat2022bregman}.
For a boundary point $\overline\z$, define the affine manifold
$
    \cM(\overline\z)
    :=
    \{\y\in\R^n:\y_{\cB(\overline\z)}
    =
    \overline\z_{\cB(\overline\z)}\},
$
where $\cB(\overline\z)$ is the set of boundary coordinates of $\overline\z$. 
Theorem~\ref{th:Q} shows that if interior feasible points approach $\overline\z$, then a vanishing $R_\gamma^t$ certifies stationarity only on $\cM(\overline\z)$, not stationarity for~\eqref{eq:obj}. 
In other words, the limiting Bregman geometry freezes the boundary coordinates. 
This mismatch gives rise to \emph{spurious stationary points}: Feasible points that are stationary on $\cM(\overline\z)$ but not stationary for~\eqref{eq:obj}; see Definition~\ref{def:spurious}. 
When $h$ has full domain, $\cM(\overline\z)=\R^n$, so the manifold and original stationarity notions coincide.


This variational mismatch has direct algorithmic consequences. 
Near a spurious stationary point, the Bregman update becomes asymptotically close to the identity map, so the iterates may move only negligibly even though the point is not stationary for the original problem.
In Theorem~\ref{th:hard}, we prove that every spurious stationary
point admits arbitrarily long local trapping: For any prescribed horizon
$K\in\N$ and tolerance $\epsilon>0$, there exists an interior feasible
initialization sufficiently close to the spurious stationary point such that
all iterates up to time $K$ remain within its $\epsilon$-neighborhood.
The result shows that the trapping phenomenon is not specific to the LP
construction, but follows from the near-identity behavior of the Bregman update
near spurious stationary points.

Theorem~\ref{th:hard} shows that trapping occurs once the initialization is sufficiently close to a spurious stationary point. 
This leaves open whether such behavior is merely a consequence of choosing the initialization artificially close to the spurious stationary point.
In Section~\ref{sec:example}, we show that this need not be
the case. We construct a nonconvex instance for which BPG, initialized at a
well-behaved interior point far from any spurious stationary point, is
nevertheless drawn toward a spurious stationary point and remains near it for
an arbitrarily long finite time. On the same instance, Euclidean projected
gradient dynamics quickly moves toward the global minimizer. Thus, the
phenomenon reflects a genuine difference between Euclidean and Bregman
dynamics, rather than only a poor choice of initialization.

It remains to ask whether the spurious points underlying these trapping
phenomena are merely artifacts of special constructions. Theorem~\ref{th:LP}
shows that they are not: For polyhedral feasible sets in standard form and
kernels whose closed domains coincide with the nonnegative orthant, every
non-stationary vertex is a spurious stationary point. Hence, spurious
stationarity is not an isolated pathology, but a structural boundary phenomenon
for Bregman geometry over polyhedral feasible sets.


Overall, our results reveal a structural blind spot in Bregman-residual-based convergence analysis for BPs with boundary-divergent kernels. 
Near the boundary, small Bregman residuals may certify only restricted stationarity, while the iterates can remain near non-stationary points for an arbitrarily long finite time. 
This motivates a more careful interpretation of Bregman stationarity measures and the development of alternative certificates or algorithmic safeguards.


	\vspace{2mm}
		



\paragraph{Notation.}  
We denote by $\overline{\R}, \, \R, \, \R_+$ the sets of extended real numbers, real numbers, and nonnegative real numbers, respectively. 
For a vector $\x \in \R^n$, its $i$-th coordinate is represented by $x_i$, and $\x_\mathcal{I}$ denotes a subvector of $\x$ indexed by $\mathcal{I}$.
The Euclidean ball $\B_{\epsilon}(\x)$ is defined as $\B_{\epsilon}(\x) := \{\y \in \R^n : \|\x - \y\| \leq \epsilon\}$.
Given a set $\mathcal{X} \subseteq \R^n$, we use $\cl(\mathcal{X})$, $\intm(\mathcal{X})$, and $\bd(\mathcal{X})$ to denote its closure, interior, and boundary, respectively.
The indicator function $\delta_{\mathcal{X}}$ of a set $\mathcal{X}$ is defined as $\delta_{\mathcal{X}}(\x) = 0$ if $\x \in \mathcal{X}$; $\delta_{\mathcal{X}}(\x) = +\infty$ otherwise. 
Unless otherwise specified, the sequence $\{\x^k\}_{k \in\N}$ always refers to the iterates generated by BPs. 
In contrast, the sequences $\{\y^k\}_{k \in\N}$ and $\{\z^k\}_{k \in \N}$ are auxiliary sequences introduced solely for technical analysis.

\paragraph{Organization.}
The remainder of the paper is organized as follows. Section~\ref{sec:pre} introduces the problem setup and necessary preliminaries. In Section~\ref{sec:main}, we present our three main theoretical results. Section~\ref{sec:example} 
provides two illustrative examples that demonstrate the implications of finite-time trapping and spurious stationarity.
The proofs of our main results are given in Section~\ref{sec:proof}. We conclude with final remarks in Section~\ref{sec:conclusion}.

\section{Preliminaries and Problem Setup }\label{sec:pre}
In this section, we introduce the key assumptions and definitions that form the foundation of our analysis. We begin with the definition of a separable kernel function.
\begin{definition}\label{def:kernel}
	A function $h:\R^n \rightarrow \overline\R$ is called a separable kernel function if it satisfies the following conditions: 
	\begin{enumerate}[label={{\rm (\roman*)}}, itemsep=0.5pt]
		\item There exists a univariate function  $\varphi:\R\rightarrow  \overline\R$ such that 
		$h(\x)=\sum_{i=1}^n\varphi(x_i),$
		where $\varphi$ is continuously differentiable on ${\rm{int}}(\dom(\varphi))$.
		\item For every sequence  $\{x^k\}_{k \in \N} \subset \mathrm{int}(\dom(\varphi)) $  converging to a point $x \in \mathrm{bd}(\dom(\varphi))$, we have $|\varphi'(x^k)| \to +\infty$. 
		\item The function $\varphi$ is strictly convex.
	\end{enumerate}
\end{definition}
The separability structure in property~{\rm (i)} is prevalent in a broad range of applications; see, e.g., \citep{bauschke2019linear, azizian2022rate, liconvergent}. Properties~{\rm (ii)} and~{\rm (iii)} are collectively referred to as \emph{Legendre-type conditions}, as introduced in \citet[Chapter 26]{rockafellar1970convex}. The following are common examples of  kernel functions that satisfy Definition~\ref{def:kernel}:
\begin{example}(See \citep[Example 1]{bauschke2017descent}.)
	\begin{enumerate}[label={{\rm (\roman*)}}, itemsep=0.5pt]
		\label{example:h}
		\item  Boltzmann–Shannon entropy kernel 
		$h(\x)=\sum_{i=1}^nx_i\log(x_i)$;
		\item Fermi–Dirac entropy kernel
		$h(\x)=\sum_{i=1}^nx_i\log(x_i)+(1-x_i)\log(1-x_i)$;
		\item Burg entropy kernel
		$h(\x)=\sum_{i=1}^n-\log(x_i)$;
		\item Fractional power kernel $h(\x)=\sum_{i=1}^npx_i-\frac{x_i^p}{1-p}$ ($0<p<1$);
		\item Hellinger entropy kernel $h(\x)=\sum_{i=1}^n-\sqrt{1-x_i^2}$.
	\end{enumerate}
\end{example}
\begin{remark}\label{remark:separable}
	The separable structure of $h$ implies that  $\cl(\dom(h))$ is a box of the form 
	\[\cl(\dom(h))=[a,c]\times[a,c]\times\cdots\times[a,c],\] 
	where $a,c\in{\R}\cup\{\pm\infty\}$, and $\cl(\dom(\varphi))=[a,c]$. 
	By the strict convexity of $\varphi$, its derivative $\varphi^\prime$ is strictly increasing. Moreover, by Property (ii) in Definition~\ref{def:kernel}, $\varphi^\prime$ diverges to $\pm \infty$ near the boundary of $\dom(\varphi)$. Throughout, we adopt the convention that  $\varphi^{\prime}(a)=-\infty$ when $a>-\infty$ and $\varphi^{\prime}(c)=+\infty$ when $c<+\infty$. 
\end{remark}

While Definition~\ref{def:kernel} imposes structural assumptions on the kernel function $h$, it also implies a continuity-like property that is frequently used in convergence analysis. Unlike many previous works that directly assume this as a technical condition, see, e.g., \citep[Definition 2.1 (vi)]{de1986relaxed}, \citep[Definition 2.1 (v)]{chen1993convergence}, \citep[Assumption H(ii)]{bauschke2017descent}, \citep[B5, p. 95]{byrne2001proximity}, \citep[Definition 3.2 (B4)]{souza2010proximal}, and \citep[Assumption H (iii)]{teboulle2018simplified}, we derive it here as a consequence of the separability and strict convexity properties of $h$. 

\begin{lemma}
	\label{lemma:bregman_limit}
	Let $h:\R^n \rightarrow \overline{\R}$ be a separable kernel function. Suppose sequences $\{\y^k\}_{k\in\N},\{\z^k\}_{k\in\N}\subseteq\intm(\dom(h))$ satisfy $\z^k\to\overline{\z}$ and $D_h(\y^k,\z^k)\to0$.  Then it follows that $\y^k\to\overline{\z}$.
\end{lemma}
\begin{proof}
	We prove the result by contradiction. Suppose, on the contrary, that $\y^k \not\to \overline{\z}$. Then, by passing to a subsequence if necessary, we may assume $y^k_{i_0}\to \overline{y}_{i_0}\in\cl(\dom(\varphi))$ with $\overline{y}_{i_0} \neq \overline{z}_{i_0}$ for some $i_0\in[n]$. Without loss of generality (WLOG), we may assume $\overline{z}_{i_0}<\overline{y}_{i_0}$. Then, there exist scalars $p,q \in \intm(\dom(\varphi))$ such that $\overline{z}_{i_0}<p<q<\overline{y}_{i_0}$. By the convergence $z_{i_0}^k \to \overline{z}_{i_0}$ and $y_{i_0}^k \to \overline{y}_{i_0}$, there exists $k_0> 0$ such that for all $k \geq k_0$, we have $z_{i_0}^k < p < q < y_{i_0}^k.$
	
	By the three-point identity for Bregman divergences (see \cite[Lemma 3.1]{chen1993convergence}), for any $x,y,z\in \dom(\varphi)$ with $z\leq x\leq y$, we have 
	\[ D_{\varphi}(z,x)+D_{\varphi}(x,y)-D_{\varphi}(z,y)=(z-x)(\varphi^{\prime}(y)-\varphi^{\prime}(x))\leq0, \]
	where the inequality follows from the strict convexity of $\varphi$. As a result, we obtain 
	\begin{equation}
		\label{eq:three-point-identity}
		D_{\varphi}(z,x)\leq D_{\varphi}(z,y) \text{ and }   D_{\varphi}(x,y)\leq D_{\varphi}(z,y), \quad \text{ for } z\leq x\leq y.
	\end{equation}
	It follows that $D_h(\y^k,\z^k)\geq D_{\varphi}(y^k_{i_0},z^k_{i_0})\geq D_{\varphi}(q,z^k_{i_0})\geq D_{\varphi}(q,p)>0$, which contradicts $D_h(\y^k,\z^k)\to0$. We complete the proof. 
\end{proof}

We now state the assumptions imposed on the optimization problem \eqref{eq:obj}. 
\begin{assumption}
	\label{assum:h} Let $\dom(g)=\cX$ be a nonempty closed convex set. We make the following assumptions:
	\begin{enumerate}[label={{\rm (\roman*)}}, itemsep=0.5pt]
		\item The function $f$ is continuously differentiable on $\cX$.
		\item The function $ g$ is convex and locally Lipschitz continuous on $\cX$. 
		\item There exists a strictly feasible point $\x^{{\rm{int}}}\in {\rm{int}}(\dom(h))\cap \cX$, and $\cX\subseteq\cl(\dom(h))$.
		\item The function $h$ is a separable kernel function; see Definition \ref{def:kernel}. 
	\end{enumerate}
\end{assumption}
Assumptions~\ref{assum:h}~{\rm(i)}–{\rm(iii)}, or their stronger variants, are standard in the literature; see, e.g.,   \citep[Assumption~A]{bauschke2017descent, bauschke2019linear}, and  \citep[Definition~1 and Assumption~1]{azizian2022rate}.

Next, we state the assumptions on the surrogate model $\gamma$ used in the update rule \eqref{eq:unified} of BPs.


\begin{assumption}[Surrogate model $\gamma$]\label{assum:gamma} 
	The following conditions hold:
	\begin{enumerate}[label={{\rm (\roman*)}}, itemsep=0.5pt]
		\item The mapping $(\x,\y)\mapsto\gamma(\y;\x)$, as well as the gradient mapping $(\x,\y)\mapsto\nabla\gamma(\y;\x)$ are jointly continuous with respect to 
		$(\y,\x)$ for all $\y\in\cX$ and $\x\in\cX$.
		\item For all $\x\in\cX$, we have $\nabla \gamma(\y;\x)\mid_{\y=\x}=\nabla f(\x), \text{ and  } \gamma(\y;\x)\mid_{\y=\x}= f(\x).$
		\item 
		There exists a constant $\overline{t}>0$ such that, for all $\x\in\cX$,  the function $\overline{t}(\gamma(\ \cdot\ ;\x)+g(\cdot))+h(\cdot)$ is strictly convex.   
		\item  
		Either $\cX$ is compact or the following condition holds: For all step sizes $t\in(0,\overline{t}]$ and all sequences $\{\z^k\}_{k\in\N},\{\y^k\}_{k\in\N}\subseteq\intm(\dom(h))\cap\cX$ with $\|\y^k\|\to+\infty$ and $\z^k\to\overline{\z}\in\cX$, we have 
		\begin{equation}
			\label{eq:assum_gamma}\lim\limits_{k\to\infty}\gamma(\y^k;\z^k)+g(\y^k)+\frac1t D_h(\y^k,\z^k)=+\infty. 
		\end{equation}
	\end{enumerate}
	Unless otherwise specified, the step size $t$ in this paper is assumed to satisfy $t\in(0,\overline{t})$.
\end{assumption} 

Assumptions~\ref{assum:gamma}~(i) and~(ii) are standard, serving to ensure the continuity and local accuracy of the surrogate model $\gamma$. In all three choices of $\gamma$ discussed in the introduction, Assumption~\ref{assum:gamma} (iii) is either a standard condition in the literature or is automatically satisfied. Specifically, when  $\gamma$ is the first-order expansion of $f$ at the current iterate $\x$, Assumption~\ref{assum:gamma} (iii) holds trivially. When $\gamma$ is taken as the original function $f$, the condition reduces to the relative convexity, a weaker assumption that has been extensively studied; see, e.g., \citep{bolte2018first,zhang2018convergence}. In the case where $\gamma$ is the second-order expansion of $f$, the $L$-smoothness of $f$ and the strong convexity of $h$ together suffice to guarantee Assumption~\ref{assum:gamma}~(iii).
Assumption~\ref{assum:gamma}~(iv) ensures the well-posedness of the BPs. If $\cX$ is compact, this condition holds automatically;   see, e.g.,  \cite[Lemma 2]{bauschke2017descent} and \citep[Assumption B]{bolte2018first}. Interested readers are referred to Appendix~\ref{appen:verify} for the rigorous verification of this assumption for commonly used surrogate models.



Following \cite[Lemma 2.3]{teboulle2018simplified}, we are now ready to state a key result concerning the well-posedness of the update rule~\eqref{eq:unified} over $\intm(\dom(h))$. 

\begin{lemma}\label{le:well}
	Suppose that Assumptions~\ref{assum:h} and~\ref{assum:gamma} hold. Then for all $\x \in \intm(\dom(h)) \cap \cX$, the update mapping $T_{\gamma}^t:\R^n\rightarrow \R^n$ defined by
	\[
	T_{\gamma}^t(\x) := \argmin_{\y \in \cX} \left\{ \gamma(\y; \x) + g(\y) + \frac{1}{t} D_h(\y, \x) \right\}
	\]
	is well-defined, and satisfies
	$
	T_{\gamma}^t(\x) \in \intm(\dom(h)) \cap \cX.
	$
\end{lemma}

\begin{proof}
	By Assumption~\ref{assum:gamma} (iv), for any sequence $\{ \y^k \}_{k\in\N}  \subset \intm(\dom(h)) \cap \mathcal{X} $ with $ \| \y^k \| \to \infty $ and $\x\in\cX$, it holds that
	$$
	\lim_{k \to \infty} \left\{ \gamma(\y^k; \x) + g(\y^k) + \frac{1}{t} D_h(\y^k, \x) \right\} = +\infty.
	$$
	This implies that the objective function is coercive over $ \intm(\dom(h)) \cap \mathcal{X} $, i.e., it tends to infinity as $ \| \y \| \to \infty $. Hence, there exists a radius $ r > 0 $ such that the infimum is attained over the compact set $ \mathcal{X} \cap \B_r(\bz) $. That is,
	$$
	\inf_{\y \in \mathcal{X}} \left\{ \gamma(\y; \x) + g(\y) + \frac{1}{t} D_h(\y, \x) \right\}
	=
	\inf_{\y \in \mathcal{X} \cap \B_r(\bz)} \left\{ \gamma(\y; \x) + g(\y) + \frac{1}{t} D_h(\y, \x) \right\}
	> -\infty.
	$$
	Since the objective function is continuous over $\cX$ and the feasible set is closed and bounded, the infimum is attained. By the definition of $ T_\gamma^t $ and \citet[Lemma 2.3 (a)]{teboulle2018simplified}, we conclude that $ T_\gamma^t(\x) $ is well-defined and lies in $ \operatorname{int}(\operatorname{dom}(h)) \cap \mathcal{X} $.
\end{proof}


Finally, to evaluate how close a given iterate is to stationarity, it is common to introduce a residual mapping $R : \mathbb{R}^n \to \mathbb{R}_+$ that quantifies the degree of stationarity.  Such a residual mapping plays a central role not only in convergence analysis but also in the design of practical stopping criteria. In this paper, we adopt the following Bregman divergence-based stationarity measure:
\begin{equation}
	R_\gamma^t(\x) := D_h\left(T^t_\gamma(\x),\x\right),
	\label{eq:stationary}
\end{equation}
which measures the discrepancy between the current iterate $\x$ and its update $T^t_\gamma(\x)$ under the Bregman geometry induced by $h$. This formulation unifies various residual-type stationarity measures commonly used in the analysis of BPs; see, e.g., \citep{bedi2022hidden, huangbregman, huang2022enhanced, latafat2022bregman}.
In particular, if we set $\gamma = f$, then the update mapping $T^t_\gamma(\x)$ reduces to the standard Bregman proximal operator, and $R_\gamma^t(\x)$ coincides with the stationarity gap $D_h(\prox_{h,F}^t(\x), \x)$ introduced by \citet{zhang2018convergence}. To make the connection precise, we recall the definition of the Bregman proximal mapping:
\begin{definition}[Bregman proximal mapping~\citep{bauschke2018regularizing,lau2022bregman}]
	Let $F : \intm(\dom(h)) \to \R$ and $t > 0$. The Bregman proximal mapping for $F$ associated with the kernel $h$ is defined by
	\[
	{\prox}^{t }_{h,F}(\x)=\underset{\y \in \R^n}{\argmin} \left\{F(\y)+\frac{1}{t} D_h(\y, \x)\right\}.
	\]
\end{definition}
Note that the function $\gamma$ in \eqref{eq:stationary} is independent of the algorithmic update and can be chosen differently. For example, as in \citet{zhang2018convergence}, in the analysis of BPG, one may let 
$\gamma = f$ in \eqref{eq:stationary} even if the algorithm uses a linear approximation of $f$ as the surrogate function instead.


\section{Main Results}\label{sec:main}
In this section, we present our main theoretical results. We begin by providing a complete characterization of the stationarity measure introduced in~\eqref{eq:stationary}. To state our main results precisely, we first define the following index sets associated with a point $\x \in \R^n$:
\[\cB(\x):=\left\{b\in [n]:x_b\in{\rm{bd}}\left(\dom(\varphi)\right)\right\}, \quad\cI(\x):=\left\{i \in [n]:x_i\in{\rm{int}}\left(\dom(\varphi)\right)\right\}.\] 
Here, 
$\cB(\x)$ and $\cI(\x)$ denote the sets of coordinate indices where 
$\x$ lies on the boundary and in the interior of $\dom(\varphi)$, respectively. We associate with $\x$ the affine manifold
\[\cM(\x)=\left\{\y\in \R^n :\y_{\cB(\x)}=\x_{\cB(\x)}\right\}.\] For any $\y\in\mathcal M(\x)$, its tangent space is
$
T_{\mathcal M(\x)}(\y)=\{\v\in\mathbb R^n:\v_{\mathcal B(\x)}=\bz\}.
$
Since $\mathcal M(\x)$ is affine, this tangent space is independent of $\y$, and we write it simply as
$T_{\cM(\x)}$.

For a base point $\x\in\cl(\dom h)$, we define the subdifferential of $g$ restricted to $\mathcal M(\x)$ by
\[
\partial_{\mathcal M(\x)}g(\u)
=
\left\{
\v\in T_{\mathcal M(\x)}:
g(\y)-g(\u)\ge \v^\top(\y-\u),\ \forall \y\in\mathcal M(\x)
\right\},
\qquad \u\in\mathcal M(\x).
\]
Similarly, we define the gradient of $f$ restricted to $\cM(\x)$ by
\begin{equation}\label{eq:manifold_gradient}
    \nabla_{\cM(\x)} f(\u)
    :=
    \proj_{T_{\cM(\x)}}(\nabla f(\u)),
    \qquad \u\in \cM(\x).
\end{equation}
These definitions coincide with the standard Riemannian subdifferential and
gradient on the affine manifold $\cM(\x)$; see
\cite[Definition~2.3]{ferreira1998subgradient} and
\cite[Proposition~3.61]{boumal2023introduction}.

The following theorem gives an exact characterization of the limiting stationarity notion certified by $R_\gamma^t$ at boundary points of the kernel domain.

 
\begin{theorem}\label{th:Q}
The following statements are equivalent for $\overline{\z} \in \cX$:
\begin{enumerate}[(i)] 
\item For every sequence $\{\z^k\}_{k\in\mathbb{N}} \subset \operatorname{int}(\operatorname{dom} h) \cap \cX$ 
with $\z^k \to \overline{\z}$, we have 
$
\lim_{k\to\infty} R^t_{\gamma}(\z^k) = 0.
$ 
\item  $\mathbf{0} \in \nabla_{\cM(\overline{\z})} f(\overline{\z}) + \partial_{\cM(\overline{\z})} g(\overline{\z})$.
\end{enumerate}
\end{theorem}
\begin{remark}[Uniform versus pathwise residual vanishing]\label{rem:uniform_pathwise}
The quantifier ``for every sequence'' in Theorem~\ref{th:Q} is essential. 
The theorem is a local variational characterization of the residual near 
$\overline{\z}$, not a pathwise convergence statement for a particular algorithmic trajectory.
This uniform formulation is also what allows us to derive the finite-time 
trapping result in Theorem~\ref{th:hard}: Near a spurious stationary point, it 
implies that the Bregman update map becomes asymptotically fixed. 
In contrast, the existence of a single interior feasible sequence 
$\z^k\to\overline{\z}$ with $R_\gamma^t(\z^k)\to0$ does not, in general, 
imply the manifold stationarity condition in Theorem~\ref{th:Q}(ii). 
Appendix~\ref{sec:tight_th:Q} gives a counterexample.
\end{remark}

Theorem~\ref{th:Q} shows that the zero-limit behavior of $R_\gamma^t$ near the boundary certifies stationarity on the affine manifold $\cM(\overline{\z})$, rather than stationarity for the original problem. 
Thus, a point may appear stationary to the Bregman update because the boundary coordinates are effectively frozen, even though it is not stationary for \eqref{eq:obj}. 
Consequently, the smallness of $R_\gamma^t$ alone cannot be interpreted as approximate stationarity for the original problem.

This limitation is not merely an artifact of the residual. 
The measure $R_\gamma^t$ records the actual Bregman displacement generated by the update and is closely tied to the descent mechanism used in standard analyses of BPs; see, e.g., \citep{bauschke2017descent,zhang2018convergence}. 
Theorem~\ref{th:Q} therefore reveals a mismatch between the local geometry measured by Bregman displacements near the boundary and the first-order stationarity condition of \eqref{eq:obj}.

We formalize this mismatch through the notion of \emph{spurious stationary points}, which arise naturally from the equivalence established in Theorem~\ref{th:Q}.



\begin{definition}[Spurious stationary points]\label{def:spurious}
	A point $\x \in \cX$ is defined as a \textit{spurious stationary point} of problem \eqref{eq:obj} if $\bz\in \nabla_{\cM(\x)}f(\x)+\partial_{\cM(\x)}g(\x)$ but $\bz \notin \partial F(\x)$.
\end{definition}

Clearly, a necessary condition for spurious stationarity is that
$T_{\cM(\x)}$ be a strict subspace of $\mathbb R^n$, equivalently
$\cB(\x)\neq\emptyset$. Thus, spurious stationary points can arise only at
boundary points of the kernel domain. When $h$ has full domain,
$\cM(\x)=\mathbb R^n$ for all $\x$, so the manifold and Euclidean
stationarity notions coincide. Moreover, if $h$ has Lipschitz continuous gradient, then by Definition~\ref{def:kernel} (ii), we have  $\dom(h)=\R^n$ and $\cI(\x)=[n]$ hold for all $\x\in\cX$, so the manifold $\cM(\x)$ coincides with the ambient space. In this sense, the Riemannian and Euclidean subdifferentials agree, and spurious stationary points cannot arise by Definition~\ref{def:spurious}.


We next present two sufficient conditions for identifying spurious stationary points, both of which follow naturally from Definition~\ref{def:spurious}.  
\begin{proposition}[Sufficient conditions for spurious stationarity]\label{prop:sufficient-spurious}
	Let $\x \in \cX$ be a non-stationary point of $F$. Then $\x$ is a spurious stationary point if either of the following conditions holds:
	\begin{enumerate}[(i)]
		\item There exists a vector $\p \in \partial F(\x)$ satisfying $\p_{\cI(\x)} = 0$;
		\item The manifold $\cM(\x)$ intersects $\cX$ only at $\x$, i.e., $\cX \cap \cM(\x) = \{\x\}$.
	\end{enumerate}
	Furthermore, if $g$ is an indicator function of a polyhedron, i.e., $g = \delta_{\cX}$ with $\cX = \{ \x \in \mathbb{R}^n : \bA\x \le \bm{b} \} \ne \emptyset$ and $\bA \in \mathbb{R}^{m \times n}$, then condition~(i) is also necessary for $\x$ to be a spurious stationary point.
\end{proposition}
\begin{remark}
	Conditions~(i) and (ii) in Proposition~\ref{prop:sufficient-spurious} capture two distinct mechanisms through which spurious stationary points can arise.
	Condition~(i) captures cases where the subgradient has nonzero components only in directions that are "invisible" to the manifold, i.e., directions in which the manifold does not permit movement. As a result, the projection of the subgradient onto the tangent space is zero, making the manifold stationarity condition appear satisfied—even though the full Euclidean stationarity fails.  Condition~(ii) corresponds to a geometric degeneracy: If the manifold $\mathcal{M}(x)$ intersects the feasible region $\cX$ only at the point $x$, then the manifold-based subdifferential condition becomes vacuously true. 
\end{remark}

\begin{proof}
(i) Since $T_{\cM(\x)}=\{\v\in\R^n:\v_{\cB(\x)}=\bz\}$,
the condition $\p_{\cI(\x)}=\bz$ gives 
$
    \proj_{T_{\cM(\x)}}(\p)=\bz .
$
Let $\p=\nabla f(\x)+\v$ with $\v\in\partial g(\x)$. For any $\y\in\cM(\x)$, the defining relation of $\cM(\x)$ gives
$\y_{\cB(\x)}=\x_{\cB(\x)}$. Hence
$
    (\y-\x)_{\cB(\x)}=\bz,
$
and therefore $\y-\x\in T_{\cM(\x)}$. Therefore,
\[
    g(\y)-g(\x)
    \ge \v^\top(\y-\x)
    =
    \proj_{T_{\cM(\x)}}(\v)^\top(\y-\x),
\]
which implies 
$
    \proj_{T_{\cM(\x)}}(\v)\in\partial_{\cM(\x)}g(\x).
$
Together with 
$
    \nabla_{\cM(\x)}f(\x)
    =
    \proj_{T_{\cM(\x)}}(\nabla f(\x)),
$
we get
$
    \bz \in 
\nabla_{\cM(\x)}f(\x)+\partial_{\cM(\x)}g(\x).
$

	
	(ii) Suppose that $\cX\cap\cM(\x)=\{\x\}$. Since $\dom(g)=\cX$, the restriction
of $g$ to $\cM(\x)$ is finite only at $\x$. Hence, for every 
$\v\in T_{\cM(\x)}$,
\[
    g(\y)-g(\x)\ge \v^\top(\y-\x),\qquad \forall \y\in\cM(\x).
\]
Indeed, both sides are zero at $\y=\x$, whereas $g(\y)=+\infty$ for every
$\y\in\cM(\x)\setminus\{\x\}$. By the definition of
$\partial_{\cM(\x)}g(\x)$, this gives
$
    \partial_{\cM(\x)}g(\x)=T_{\cM(\x)}.
$
Therefore,
\[
    \bz\in \nabla_{\cM(\x)}f(\x)+\partial_{\cM(\x)}g(\x).
\]

(iii) We now prove the necessity of  Condition~(i) when $g=\delta_{\cX}$ and $\cX$ is polyhedral. 
By the definition of $\partial_{\cM(\x)}g(\x)$ and the fact that
$\y-\x\in T_{\cM(\x)}$ for all $\y\in\cM(\x)$, we have 
	\begin{equation}\label{eq:Mg}
		\begin{aligned}
			\partial_{\cM(\x)}g(\x)&=\left\{\proj_{T_{\cM(\x)}}(\v):g(\y)-g(\x)\geq\v^{\top}(\y-\x),~\forall~\y\in\cM(\x)\right\}\\
			&=\left\{\proj_{T_{\cM(\x)}}(\v):\v\in\partial (g+\delta_{\cM(\x)})(\x)\right\}.
		\end{aligned}
	\end{equation} 
	Since $g=\delta_{\cX}$, it follows that $g+\delta_{\cM(\x)} = \delta_{\cX\cap \cM(\x)}$. Note that  
	\[
	\cX\cap\cM(\x)=\{\y\in\R^n:\bA\y\leq\b,~\y_{\cB(\x)}=\x_{\cB(\x)}\}, 
	\]
	is a polyhedron, hence the subdifferential of $g+\delta_{\cM(\x)}$ at $\x$ is given by 
	\[\partial (g+\delta_{\cM(\x)})(\x)=\left\{\bA^{\top}\blam+\v:\blam\in\R^m_+,\blam^{\top}(\bA\x-\b)=0 ,\v_{\cI(\x)} = \bz\right\}. \]
Substituting this expression into \eqref{eq:Mg} and using $T_{\cM(\x)}(\x)=\{\v\in\R^n:\v_{\cB(\x)}=\bz\}$, we obtain
	\[\begin{aligned} 
		\partial_{\cM(\x)}g(\x)&=\left\{\proj_{T_{\cM(\x)}}(\bA^{\top}\blam):\blam\in\R^m_+,\blam^{\top}(\bA\x-\b)=0\right\}.
	\end{aligned}\]
   It leads to
	\[
	\bz\in \nabla_{\cM(\x)}f(\x)+\partial_{\cM(\x)}g(\x) \, 
	\Longleftrightarrow \,  \bz\in \left\{\proj_{T_{\cM(\x)}}(\nabla f(\x)+\bA^{\top}\blam):\blam\in\R^m_+,\blam^{\top}(\bA\x-\b)=0 \right\}. 
	\]
	Equivalently, 
	\[
	\bz\in \left\{\proj_{T_{\cM(\x)}}(\p):\p\in\partial F(\x)\right\},
	\]
	which is precisely Condition~(i). We complete the proof. 
\end{proof}

While the previous discussion exposes a limitation of Bregman-residual-based
stationarity certificates, the following result shows that the issue is also
reflected in the update dynamics. Once a spurious stationary point exists, the
Bregman update can remain arbitrarily close to it for any prescribed finite
number of iterations.


\begin{theorem}[Finite-time trapping near spurious stationary points]\label{th:hard}
	Suppose that there exists a spurious stationary point $\tilde{\x}^\star \in \cX$ for problem \eqref{eq:obj}. For every $K\in\N$ and $\epsilon>0$, there exists an initial point $\x^0\in \B_{\epsilon}(\tilde{\x}^\star)\cap \cX\cap\intm(\dom(h))$, sufficiently close to the spurious stationary point $\tilde{\x}^\star$, such that the sequence $\{\x^k\}_{k\in[K]}$ generated by \eqref{eq:unified} with a fixed step size $t\in(0,\overline{t})$ satisfies
	\begin{equation}\label{eq:trap}
		\x^k\in\B_{\epsilon}(\tilde{\x}^\star)\quad\text{ for all }k\in [K].
	\end{equation} 
\end{theorem}


\begin{proof}
	By Theorem~\ref{th:Q} and Definition \ref{def:spurious}, for every sequence $\{\z^k\}_{k\in\N}\subseteq\intm(\dom(h))\cap\cX$ converging to a spurious stationary point $\tilde{\x}^\star$, we have $\lim_{k\to\infty} R_\gamma^t(\z^k) = \lim_{k\to\infty} D_h(T^t_{\gamma}(\z^k),\z^k) =0$. 
	Moreover, under Assumption \ref{assum:h} (iv), Lemma \ref{lemma:bregman_limit} implies that $\lim_{k\to\infty}{T}^t_{\gamma}(\z^k)=\tilde{\x}^\star$. 
Hence, for any $r>0$, there exists $\eta(r)>0$ such that for all 
	$\x\in\cX\cap\intm(\dom(h))$ with 
	$\|\x-\tilde{\x}^\star\|<\eta$, we have
	$	    \|T^t_{\gamma}(\x)-\tilde{\x}^\star\|<r .
	$
	Set $\epsilon_0:=\epsilon$. For each $k=0,\ldots,K-1$, applying the above
local property with $r=\epsilon_k$, choose
$
    0<\epsilon_{k+1}
    <
    \min\left\{\eta(\epsilon_k),\frac{\epsilon_k}{2}\right\}.
$
Then, we have
\[
    \|\x-\tilde{\x}^\star\|<\epsilon_{k+1}
    \quad\Longrightarrow\quad
    \|T_\gamma^t(\x)-\tilde{\x}^\star\|<\epsilon_k .
\]
Thus, if the initial point $\x^0\in\cX\cap\intm(\dom(h))$ satisfies $
	    \|\x^0-\tilde{\x}^\star\|<\epsilon_K,
	$ 
	then induction gives
	\[
	    \|\x^j-\tilde{\x}^\star\|<\epsilon_{K-j}
	    \le \epsilon_0=\epsilon,
	    \qquad j=0,1,\ldots,K.
	\]
	Such an initial point exists by Assumption~\ref{assum:h}(iii) and the convexity of $\cX$.
	Hence all iterates up to time $K$ remain in $\B_\epsilon(\tilde{\x}^\star)$. 
	This completes the proof.
	

\end{proof}
Theorem~\ref{th:hard} reduces the finite-time trapping phenomenon to the local
presence of spurious stationary points. It is therefore important to understand
whether such points are merely artifacts of special constructions. The next
result shows that they are not: For constrained nonconvex problems over
standard-form polyhedral feasible sets, and for kernels whose closed domains
coincide with the nonnegative orthant, every non-stationary vertex is a
spurious stationary point. Thus, spurious stationarity is not an isolated pathology, but a structural
boundary phenomenon arising from the interaction between polyhedral constraints
and singular Bregman geometries.
\begin{theorem}\label{th:LP}
	Consider the optimization problem 
	\[
	\min_{\x \in \mathcal{X}} ~f(\x), \quad \text{where } \mathcal{X} := \{ \x \in \R^n : \mathbf{A}\x = \mathbf{b},\ \x \geq \bz\}
	\]
	is not singleton, with $\mathbf{A} \in \R^{m \times n}$ and $\mathbf{b} \in \mathbb{R}^m$. 
    Assume that $f$ is continuously differentiable on $\cX$, and that the kernel function $h$ satisfies $\cl(\dom(h))= \R^n_+$. 
    Then, every non-stationary vertex of $\mathcal{X}$ is a spurious stationary point.
\end{theorem}

\begin{proof}
	At any vertex $\x$ of $\cX$, the active constraint gradients consist of the rows of $\mathbf{A}$ together with the standard basis vectors $\e_i$ for all $i \in \cB(\x)$. These vectors together span $\mathbb{R}^n$, that is,
	\[
	\rank\left(\left[\mathbf{A}^T,\e_i: i \in \cB(\x) \right] \right) = n,
	\]
	as shown in \cite[Section 8.5]{schrijver1998theory}.
	Since $\cX$ is not a singleton, we must have $\rank(\mathbf{A}) < n$. Hence, any vertex $\x \in \cX$ must have at least one active inequality constraint, i.e., $\mathcal{B}(\x) \neq \emptyset$.

    By the rank condition, there exist vectors $\bmu \in \mathbb{R}^m$ and $\blam \in \mathbb{R}^{|\cB(\x)|}$ such that
	\[
	 \nabla f(\x)+ \mathbf{A}^T\bmu + \sum_{i \in \cB(\x)} \lambda_i \e_i = \bz, 
	\]
	which implies that $(\nabla f(\x)+\mathbf{A}^T\bmu)_{\cI(\x)}=\bz.$ 
 Let \(F:=f+\delta_{\cX}\) and define
$
    \p:=\nabla f(\x)+\mathbf A^\top\bmu .
$
Then $\p\in\partial F(\x)$ and $\p_{\cI(\x)}=\bz$. Since $\x$ is non-stationary, $\bz\notin\partial F(\x)$. Hence, Proposition~\ref{prop:sufficient-spurious}(i) shows that $\x$ is a spurious stationary point.	
\end{proof}

\section{Practical Implications of Finite-Time Trapping} \label{sec:example}

Theorem~\ref{th:hard} shows that, near a spurious stationary point, the Bregman update can become locally almost fixed, leading to arbitrarily long finite-time trapping. 
This section illustrates two implications of this result. 
The first example revisits the LP instance in~\eqref{example:lp} and explains why this behavior does not contradict standard objective-gap guarantees for convex BPG: the relevant initial Bregman divergence can be arbitrarily large near a spurious stationary point. 
The second example shows that finite-time trapping is not merely an artifact of initializing close to a spurious stationary point; in a nonconvex problem, BPG can be drawn toward such a point from a well-behaved interior initialization. 
Together, these examples underscore the subtle yet widespread risk posed by spurious stationary points under singular Bregman geometry.



We first examine the LP instance given in \eqref{example:lp}.  

\begin{example}
	Suppose that $\cl(\dom(h))=\R^2_+$ and consider the problem 
	\[
	\begin{array}{rl}
		\min & -x_1 \\
		\text{s.t.} & x_1 + x_2 = 1, \quad x_1, x_2 \geq 0.
	\end{array}
	\]
	A simple calculation reveals that the point 
	$\tilde{\x}^\star = (0,1)$  is not a stationary point: 
	\[
	\begin{aligned}
		\bz \notin \partial F(\tilde{\x}^\star) 
		&= \left\{ (-1,0) + \lambda (-1,0) + \mu (1,1) : \lambda\in\R_+, \mu \in \mathbb{R} \right\}.
	\end{aligned}
	\]
	Moreover, the interior coordinate set at $\tilde{\x}^\star$ is $\mathcal{I}(\tilde{\x}^\star) = \{2\}$, and there exists a subgradient $\p \in \partial F(\tilde{\x}^\star)$ such that $\p_{\mathcal{I}(\tilde{\x}^\star)} = 0$. By Proposition \ref{prop:sufficient-spurious} (i), we conclude that $\tilde{\x}^\star$ is a spurious stationary point. 
\end{example}

To illustrate Theorem~\ref{th:hard}, we adopt the Boltzmann–Shannon entropy kernel $\varphi(x) = x \log x$, which is widely used in constrained optimization problems over the simplex. Under this kernel, BPG update reduces to the classical multiplicative weights update method \citep{arora2012multiplicative}:  
\begin{align*}
	\x^{k+1} & = \argmin_{\y\in\R^2} t (-1,0)^T \y + D_h(\y,\x^k) +\delta _{ \Delta_2}(\y) \\
	& = \left(\frac{x_1^k}{x_1^k+e^{-t}x_2^k}, \frac{e^{-t}x_2^k}{x_1^k+e^{-t}x_2^k}\right), \quad \forall k \in [K].
\end{align*} 
We initialize the algorithm with
\[
\x^0 = \left(\frac{\sqrt{2}\epsilon}{2} e^{-tK}, 1- \frac{\sqrt{2}\epsilon}{2} e^{-tK}\right),
\]
which lies strictly inside the simplex. 
Then, it is straightforward to verify that for all $k\in[K]$, 
\[
\|\x^{k+1}-\tilde{\x}^\star\|= \frac{\sqrt{2}x_1^k}{x_1^k+e^{-t}x_2^k} \leq \sqrt{2} e^t x_1^k \leq \sqrt{2} e^{tk} x_1^0 = e^{-t(K-k)}\epsilon \leq \epsilon,
\]
where the first inequality is derived from the constraint $x_1^k + x_2^k = 1$ and $t \geq 0$, the second inequality is justified by iteratively applying the recursive relation from the first inequality $k$ times.

This example shows that when the initial point is extremely close to a spurious stationary point, the trapping phenomenon described in Theorem~\ref{th:hard} can be triggered. Importantly, this behavior does not contradict the non-asymptotic convergence guarantee established in \citet[Corollary 1]{bauschke2017descent}.  That result states that  the sequence $\{\x^k\}_{k \in \mathbb{N}}$ generated by BPG satisfies: 
\begin{equation}\label{eq:rate}
	f(\x^K)-\min_{\x\in\cX}\ f\leq \frac{D_h(\x^\star,\x^0)}{t}\cdot\frac{1}{K},
\end{equation} 
where $\x^\star\in \argmin_{\x\in\cX}f$ is the global minimizer, $t$ is the step size, and $\x^0$ is the initial point. Here, the rate depends linearly on the initial Bregman divergence $D_h(\x^\star,\x^0)$, which can be made arbitrarily large by placing $\x^0$ near a spurious stationary point. For example, in the LP instance above, choosing $\x^0 = (\exp(-K),1-\exp(-K))$ yields $
D_h((1,0), (\exp(-K),1-\exp(-K))) = K$, leading to the trivial upper bound: 
\[f(\x^K)-\min_{\x\in\cX} f(\x) \leq \frac{1}{t}.\]
This illustrates that the standard non-asymptotic objective-gap bound can become
uninformative when the initialization approaches a spurious
stationary point: although the rate is $\cO(1/K)$ for a fixed initialization,
the initial Bregman divergence to the optimizer may scale with the prescribed
horizon $K$, making the resulting bound $O(1)$ on the constructed instance.

The LP example above may suggest that finite-time trapping only occurs when the initialization is chosen extremely close to a spurious stationary point. 
This intuition is reasonable in convex settings, where global objective-gap guarantees still govern the long-run behavior. 
For nonconvex problems, however, the situation is more subtle: the Bregman dynamics may actively drive the iterates toward a spurious stationary point from a well-behaved interior initialization. 
Once the iterates enter a sufficiently small neighborhood of such a point, the local trapping mechanism in Theorem~\ref{th:hard} can take over. 
The following example demonstrates this behavior.

\begin{example}
	\label{example:nonconvex}
	Suppose that $\cl(\dom(h))=\R^2_+$ and consider the following parameterized nonconvex optimization problem:
	\begin{equation}
		\begin{array}{rl}
			\min\limits_{\x\in\R^2}&\ f_\alpha(x_1,x_2):=\frac12\phi_\alpha(x_1)\cdot(x_2+0.05)  \\
			{\rm s.t.}& x_1,x_2\in[0,1],
		\end{array}  
		\label{eq:nonconvex_spurious}
	\end{equation}
	where $\alpha\in (0,0.1]$ and $\phi_{\alpha}:\R\rightarrow \R$ is a continuously differentiable function satisfying 
	\begin{enumerate}[label={{\rm (\roman*)}}, itemsep=0.5pt]
		\item $\phi_{\alpha}(x)=2(x-\alpha)$ for $x\geq \alpha$, $\phi_{\alpha}(x)\leq0$ for $x\leq \alpha$, and $\phi_{\alpha}(0)\leq-1$;
		\item $\phi^{\prime}_{\alpha}(x)\geq1$ for $x\in[0,1]$;
		\item  $f_{\alpha}+h$ and $-f_{\alpha}+h$ are convex on $(0,1]\times(0,1]$.
	\end{enumerate}
\end{example}
First of all, we compute the first-order optimality condition of \eqref{eq:nonconvex_spurious} as 
\begin{equation}\label{eq:partialF}
	\partial F(\x)=\left\{\left(\frac12\phi^{\prime}_{\alpha}(x_1)(x_2+0.05),\frac12\phi_{\alpha}(x_1)\right)+\boldsymbol{\lambda}-\boldsymbol{\mu}:\boldsymbol{\lambda},\boldsymbol{\mu}\geq\bz,\boldsymbol{\lambda}^{\top}(\1-\x)=0,\boldsymbol{\mu}^{\top}\x=0\right\}.
\end{equation}
Combining the subdifferential characterization in \eqref{eq:partialF} with condition~(ii) in Example \ref{example:nonconvex}, we observe that for any $\x \in (0,1] \times [0,1]$, all subgradients $\p \in \partial F(\x)$ satisfy $p_1 \geq 0.025$. This implies that any true stationary or spurious stationary point must lie on the boundary $x_1 = 0$.  Substituting $x_1 = 0$ into \eqref{eq:partialF} and invoking Proposition~\ref{prop:sufficient-spurious}~(i), we find that $\tilde{\x}^\star = (0,0)$ is a spurious stationary point, while $\x^\star = (0,1)$ is the unique true stationary point and hence the global minimizer.

Importantly, from the above discussion, it follows directly that $\dist(\boldsymbol{0}, \partial F(\x)) \geq 0.025$ for all $\x \in [0,1] \times [0,1] \setminus \{\x^\star\}$. This sharpness property \citep{burke1993weak} highlights that, in the Euclidean
stationarity geometry, every point away from the global minimizer admits a
uniform first-order stationarity gap. 
In the numerical comparison below, this property is reflected in the rapid
progress of Euclidean projected gradient dynamics toward the global minimizer. Moreover, condition~(iii) in Example~\ref{example:nonconvex} is a commonly adopted assumption to guarantee non-asymptotic convergence rates for BPG under the stationarity measure \eqref{eq:stationary}; see, e.g., \citep[Theorem 4.1]{zhang2018convergence} and \citep[Proposition 4.1]{bolte2018first}.

Given these favorable properties, the problem in Example \ref{example:nonconvex} may initially appear benign: It admits a unique global minimizer, enjoys a sharpness property that favors fast convergence under PGD, and satisfies structural conditions often used to guarantee non-asymptotic convergence for BPG with the measure \eqref{eq:stationary}. It is therefore tempting to expect that initializing sufficiently far from the
spurious stationary point should avoid the local trapping behavior.

This intuition is misleading. Despite the seemingly favorable properties discussed above, the iterates generated by BPG can still be drawn toward the spurious stationary point $\tilde{\x}^\star = (0,0)$ and become trapped. We begin by outlining the high-level intuition behind this phenomenon, followed by a concrete construction that formalizes this behavior.

We now instantiate the BPG method using the Burg entropy kernel $\varphi(x) = -\log(x)$, and derive its explicit update rule by applying the first-order optimality condition to~\eqref{eq:nonconvex_spurious}:

\begin{equation}
	\label{eq:kkt_nonconvex}
	\left\{
	\begin{aligned}
		& t \nabla f_\alpha(x_1^k, x_2^k) 
		+ \left( \frac{1}{x_1^k}, \frac{1}{x_2^k} \right) 
		- \left( \frac{1}{x_1^{k+1}}, \frac{1}{x_2^{k+1}} \right) 
		+ \boldsymbol{\lambda}^k = \boldsymbol{0}, \\[1ex]
		& \boldsymbol{\lambda}^k \geq \boldsymbol{0}, \quad 
		\left(\boldsymbol{\lambda}^k\right)^\top ( \boldsymbol{1} - \boldsymbol{x}^{k+1}) = \boldsymbol{0}.
	\end{aligned}
	\right.
\end{equation}
Since we observe that $\nabla_{x_1} f_\alpha(x_1^k, x_2^k) = \frac{1}{2} \phi_\alpha'(x_1^k)(x_2^k + 0.05) > 0$, the BPG update in the $x_1$-coordinate yields $x_1^{k+1} < x_1^k \leq 1$. This implies that the iterates are monotonically decreasing in $x_1$ and will continue contracting toward zero. The behavior along the $x_2$-coordinate, however, is more subtle. The first-order optimality condition leads to the explicit update:
\begin{equation}
	\label{eq:x2_update}
	\frac{1}{x_2^{k+1}} = \max\left(\frac{1}{x_2^k} + \frac{t}{2} \phi_\alpha(x_1^k),\ 1\right).
\end{equation}
When $\alpha$ is small, during the early iterations with $x_1^k > \alpha$, condition~(i) in Example~\ref{example:nonconvex} guarantees that $\phi_\alpha(x_1^k) > 0$. As a result, both $x_1^k$ and $x_2^k$ decrease simultaneously, driving the iterates toward the origin. However, once $x_1^k < \alpha$, the function $\phi_\alpha$ becomes nonpositive, and the $x_2$-update begins to reverse direction. Unfortunately, by this point, the iterates may have already been drawn sufficiently close to the spurious stationary point $(0,0)$, making it too late for the algorithm to recover. As a result, the sequence becomes trapped in a neighborhood of the spurious point.

\begin{figure}[htbp]
	\centering
	\includegraphics[width=1.0\textwidth]{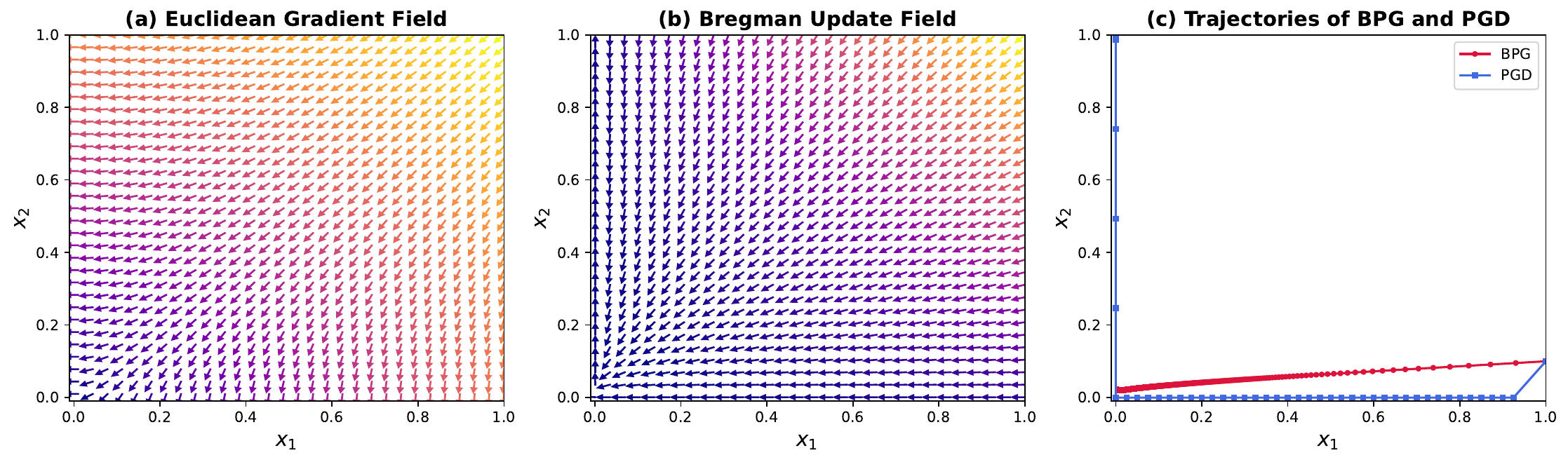}
	\caption{ Comparison of PGD and BPG dynamics with the Burg entropy kernel on the nonconvex instance constructed in Example~\ref{example:nonconvex}, where $\phi_\alpha$ is given in~\eqref{eq:concrete_phi}. 
		(a) Euclidean gradient field $-\nabla f_\alpha(\x)$. 
		(b) Bregman update field induced by Burg entropy, i.e., $T_{\gamma}^{1}(\x)-\x$. 
		(c) Trajectories of BPG and PGD starting from $(1, 0.1)$. 
		While PGD quickly converges, BPG becomes trapped near the spurious stationary point $(0,0)$.}
	\label{figure:2}
\end{figure}

We are now ready to present an explicit construction of the function $\phi_{\alpha}$ that satisfies the conditions in Example~\ref{example:nonconvex}. Specifically, we define $\phi_{\alpha}:\R \rightarrow \mathbb{R}$ as the following piecewise continuously differentiable function:
\begin{equation}
	\label{eq:concrete_phi}
	\phi_{\alpha}(x) =
	\begin{cases}
		2x - 2\alpha, & \text{if } x \in [\alpha, +\infty); \\
		x + \alpha\log\left(\frac{x}{\alpha}\right) - \alpha, & \text{if } x \in [\alpha\exp(-\frac{1}{\alpha}), \alpha]; \\
		x - \alpha\log\left(\frac{x}{\alpha}\right) + 2\exp(\frac{1}{\alpha})x - 3\alpha - 2, & \text{if } x \in [\frac{1}{2}\alpha\exp(-\frac{1}{\alpha}), \alpha\exp(-\frac{1}{\alpha})]; \\
		x + \alpha(\log 2 - 2) - 1, & \text{if } x \in (-\infty, \frac{1}{2}\alpha\exp(-\frac{1}{\alpha})].
	\end{cases}
\end{equation}
To visualize this spurious-attraction behavior in  Example~\ref{example:nonconvex}, we set $\alpha = 0.01$, use a maximum iteration count of 5000 and initialize the iterates at $(x_1^0, x_2^0) = (1, 0.1)$, which is far from both the boundary of the domain and the spurious stationary point $(0,0)$. Figure~\ref{figure:2} compares the dynamics of PGD and BPG on this instance. In Figure~\ref{figure:2} (c), we observe that PGD rapidly converges to the global minimum $(0,1)$ in fewer than 50 iterations. In stark contrast, the BPG trajectory is gradually drawn toward the spurious stationary point $(0,0)$, where it becomes trapped for more than 5000 iterations. Figure~\ref{figure:2}~(a) and~(b) illustrate the gradient fields under Euclidean and Bregman geometries, respectively. Notably, in Figure~\ref{figure:2} (b), we see that within the region $x_1 > \alpha$, BPG update directions consistently point toward $(0,0)$, making the spurious stationary point an attractor. Once the iterate enters the zone $x_1 < \alpha$, the dynamics flip abruptly due to the sign change in $\phi_{\alpha}(x_1)$, but by then the algorithm has already entered the trap. This attract-and-flip behavior illustrates how nonconvexity can drive Bregman dynamics into the trapping region associated with a spurious stationary point.

Finally, we summarize the above behavior in the following formal observation, which demonstrates that the BPG iterates can remain arbitrarily close to a spurious stationary point for any prescribed number of iterations.

\begin{observation}\label{observation:encounter2}  
	Consider Example~\ref{example:nonconvex}, where the objective function 
	$\phi_{\alpha}$  is specified as in \eqref{eq:concrete_phi}. We apply BPG with the Burg entropy kernel $\varphi(x) = -\log(x)$, using a constant stepsize $t = 0.5$ and initialization $\x^0 = (1, 0.1)$.  Let $\{\x^k\}_{k\in\N}$ denote the resulting sequence of iterates.
	Then, for any given $K \in \mathbb{N}$ and $\epsilon \in (0,1]$, there exists a parameter $\alpha > 0$ and an iteration index $k_1 > 0$ such that the iterates satisfy
	\[
	\x^k \in \mathbb{B}_\epsilon(\tilde{\x}^\star), \quad \text{for all } k_1 \leq k \leq K + k_1. 
	\]
\end{observation}
The rigorous proof of this result is provided in Appendix~\ref{appen:nonconvex_spurious}.


\section{Proof of Theorem \ref{th:Q}}\label{sec:proof}
We now present the proof of Theorem \ref{th:Q}. The overall structure follows the diagram below:
\begin{center}
	\begin{tikzpicture}[scale=0.8]
		\node[rectangle,
		minimum width =100pt ,
		minimum height =20pt ,draw=black] (1) at(-6,4){$\bz \in \nabla_{\cM(\overline{\z})}f(\overline{\z})+\partial_{\cM(\overline{\z})}g(\overline{\z})$};
		\node[rectangle,
		minimum width =20pt ,
		minimum height =20pt ,draw=black] (2) at(1,4){$\overline{ R}^t_{\gamma}(\overline{\z})=0$}; 
		\node[rectangle,
		minimum width =20pt ,
		minimum height =20pt ,draw=black] (3) at(7,4){$\lim\limits_{k\to\infty}{ R}^t_{\gamma}(\z^k)=0$};
		\draw[double, <->](1) --(2) node[midway, sloped,above] {Prop. \ref{pro:bridge}  } ;      
		\draw[double, <->] (2) --(3)  node[midway, sloped, above] {Prop. \ref{pro:R} };
	\end{tikzpicture}
\end{center}
The proof proceeds in three steps.
First, we define the extended Bregman stationarity measure $ \overline{R}_\gamma^t$, which is well-defined on the entire domain $\cX$, and satisfies $\overline{ R}^t_{\gamma}(\x)={R}^t_{\gamma}(\x)$ for $\x\in \cX\cap\intm(\dom(h))$. This measure serves as a technical bridge between the limiting behavior of the original stationarity measure   $R_\gamma^t(\z^k)$ and the variational condition that $\bz \in \nabla_{\cM(\overline{\z})}f(\overline{\z})+\partial_{\cM(\overline{\z})}g(\overline{\z})$.  

Second, we establish that $\overline{R}_\gamma^t(\overline{\z}) = 0$ if and only if $\bz \in \nabla_{\cM(\overline{\z})}f(\overline{\z})+\partial_{\cM(\overline{\z})}g(\overline{\z})$.
	\begin{proposition}\label{pro:bridge}
		For all $\overline{\z}\in\cX$, we have the following equivalence: 
		\[\overline{ R}^t_{\gamma}(\overline{\z})=0\, \Longleftrightarrow\, \overline{T}_\gamma^t(\overline{\z}) = \overline{\z} \, \Longleftrightarrow \, \bz \in \nabla_{\cM(\overline{\z})}f(\overline{\z})+\partial_{\cM(\overline{\z})}g(\overline{\z}).\]
	\end{proposition} 
The proof of Proposition~\ref{pro:bridge} is given in Section~\ref{sec:proof_pro_bridge}.

Third, we establish the following equivalence.  

\begin{proposition}\label{pro:R}
Let $\overline{\z}\in\cX$. 
Then the following are equivalent:
\begin{enumerate}[(i)]
\item $\overline{R}^t_{\gamma}(\overline{\z})=0$.
\item For every sequence $\{\z^k\}_{k\in\mathbb{N}}\subset \cX\cap\operatorname{int}(\operatorname{dom} h)$ with $\z^k\to\overline{\z}$, we have $\lim_{k\to\infty} R^t_{\gamma}(\z^k)=0$.
\end{enumerate}
\end{proposition}
	We defer the proof of Proposition~\ref{pro:R} to Section~\ref{sec:proof_conti}.

Combining Propositions \ref{pro:bridge} and \ref{pro:R} completes the proof of Theorem \ref{th:Q}. We now present the details.

\subsection{Extended Bregman Stationarity Measure}
The original stationarity measure $R_\gamma^t$ is defined only at points in
$\cX\cap\operatorname{int}(\dom h)$. Since Theorem~\ref{th:Q} concerns limits
approaching boundary points of $\dom h$, we introduce an extended Bregman
update mapping on the entire feasible set $\cX$. This extension freezes the
boundary coordinates and applies the Bregman update only along the interior
coordinates. It is used solely as an analytic device to capture the limiting
behavior of $R_\gamma^t$ near the boundary.

\begin{definition}[Extended Bregman Update Mapping]\label{def:mapping}	
	We define the extended Bregman update mapping $\overline{T}_\gamma^t(\x)$ for all $\x \in \mathcal{X}$ by
	\[
	\overline{T}_\gamma^t(\x) := \argmin_{\y \in{\R^n}} \; G_\gamma^t(\y; \x),
	\]
	where the objective $G_\gamma^t(\y; \x)$ is given by
	\[G^{t}_{\gamma}(\y;\x):= \gamma(\y;\x)+ g(\y)+\underbrace{\frac{1}{t}\sum_{i\in \cI(\x)}D_{\varphi}(y_i,x_i)}_{{\rm Interior\ coordinates}}+ \underbrace{\delta_{\cM(\x)}(\y)}_{{\rm Boundary\ coordinates}}.\]
\end{definition}
\noindent Here, the indicator term enforces the coordinates corresponding to the boundary indices $\cB(\x)$ remain fixed, while the Bregman update is applied only over the interior coordinates $\cI(\x)$.

The following result collects the basic properties of $\overline T_\gamma^t$ needed in the sequel.
\begin{proposition}\label{prop:extended-map-properties}
	The following properties hold for the extended Bregman update mapping $\overline{T}_\gamma^t$:
	\begin{enumerate}[label={{\rm (\roman*)}}, itemsep=0.5pt]
		\item \textbf{(Well-posedness)} The mapping $\overline{T}_\gamma^t(\x)$ is well-defined and single-valued for all $\x\in\cX$. 
		\item \textbf{(Boundedness)} For any bounded sequence $\{\z^k\}_{k\in\N}\subseteq\cX$, the sequence $\{\overline{T}^{t}_{\gamma}(\z^k)\}_{k\in\N}$ is also bounded.
		\item \textbf{(Boundary coordinate consistency)} Suppose that the sequence $\{\z^k\}_{k\in\N}\subseteq  \cX$ converges to ${\overline{\z}}$ and $\overline{T}^{t}_{\gamma}(\z^k)\to\overline{\v}\in\cX$. Then, we have $\overline{\v}_{\cB(\overline{\z})}=\overline{\z}_{\cB(\overline{\z})}$ and $\cB(\overline{\v})=\cB(\overline\z)$. In particular,  $\cB(\overline{T}^t_{\gamma}(\overline{\z}))=\cB(\overline\z)$.
	\end{enumerate}
\end{proposition}

The proof of Proposition~\ref{prop:extended-map-properties} is deferred to
Appendix~\ref{app:extended-update}.

Armed with the extended update mapping, we proceed to define the extended Bregman stationarity measure $\overline{ R}^t_{\gamma}$ over the entire domain $\cX$, whose well-definedness is ensured by that of $\overline{T}^t_{\gamma}$.
\begin{definition}[Extended stationarity measure]	\label{def:measure}
	We define the \textit{extended Bregman stationarity measure} $\overline{ R}^t_{\gamma}(\x):\cX \to\R_+$ as
	$\overline{ R}^t_{\gamma}(\x):= \sum_{i\in \cI(\x)}D_{\varphi}(\overline{T}^t_{\gamma}(\x)_i,x_i).$
\end{definition}

\subsection{Proof of Proposition \ref{pro:bridge}}
\label{sec:proof_pro_bridge}
\begin{proof}
	We first claim that $\overline{R}_\gamma^t(\overline{\z}) = 0 $ if and only if  $\overline{T}_\gamma^t(\overline{\z}) = \overline{\z}$.  By definition of $\overline{R}_\gamma^t$, we have:
	\[
	\overline{R}_\gamma^t(\overline{\z}) = 0 \quad \Longleftrightarrow \quad \overline{T}_\gamma^t(\overline{\z})_i = \overline{z}_i \quad \text{for all } i \in \cI(\overline{\z}).
	\]
	On the other hand, by construction of $\overline{T}_\gamma^t$, we always have
	$\overline{T}^{t}_{\gamma}(\overline{\z})_{\cB(\overline{\z})}=\overline{\z}_{\cB(\overline{\z})}$ 
	due to the hard constraint imposed on the boundary coordinates.  Therefore, $\overline{T}_\gamma^t(\overline{\z}) = \overline{\z}$ if and only if $\overline{R}_\gamma^t(\overline{\z}) = 0$.
	
	It remains to show that $ \overline{T}_\gamma^t(\overline{\z}) = \overline{\z}$ if and only if $\bz \in \nabla_{\cM(\overline{\z})}f(\overline{\z})+\partial_{\cM(\overline{\z})}g(\overline{\z})$.  From the definition of the extended Bregman update mapping,  
    we know that $\overline{T}_\gamma^t(\overline{\z}) = \overline{\z}$ if and only if $\overline{\z} \in \argmin_{\y \in \R^n} G_\gamma^t(\y; \overline{\z})$.
	By the convexity of $G_\gamma^t(\cdot; \overline{\z})$, this is equivalent to the first-order optimality condition: $\bz \in \partial G_\gamma^t(\overline{\z}; \overline{\z}).$ According to Assumption~\ref{assum:gamma} (ii) and \cite[Exercise~8.8(c)]{rockafellar2009variational}, we have
	\begin{equation}\label{eq:G_op}
		\begin{aligned}
			\partial G^t_{\gamma}(\y;\overline{\z})\mid_{\y=\overline{\z}}=\nabla f(\overline{\z})+\partial \left(g+\delta_{\cM(\overline{\z})}\right)({\overline{\z}}). 
		\end{aligned}
	\end{equation}  
We next rewrite the right-hand side of \eqref{eq:G_op} in terms of the affine
manifold $\cM(\overline{\z})$. Recall that
\[
    T_{\cM(\overline{\z})}
    =
    \{\v\in\R^n:\v_{\cB(\overline{\z})}=\bz\},
    \qquad
    T_{\cM(\overline{\z})}^{\perp}
    =
    \operatorname{span}\{\e_b:b\in\cB(\overline{\z})\}.
\]
For any $\y\in\cM(\overline{\z})$, we have
$\y-\overline{\z}\in T_{\cM(\overline{\z})}$. Hence, for any $\v\in\R^n$, we have 
$
    \v^\top(\y-\overline{\z})
    =
    \proj_{T_{\cM(\overline{\z})}}(\v)^\top(\y-\overline{\z}).
$
Therefore, the subgradient inequality defining
$\partial(g+\delta_{\cM(\overline{\z})})(\overline{\z})$ depends only on the
projection of $\v$ onto $T_{\cM(\overline{\z})}$. Indeed,
$
    \v\in\partial(g+\delta_{\cM(\overline{\z})})(\overline{\z})
$
if and only if
\[
    g(\y)-g(\overline{\z})
    \ge
    \proj_{T_{\cM(\overline{\z})}}(\v)^\top(\y-\overline{\z}),
    \qquad \forall \y\in\cM(\overline{\z}), 
\]
or equivalently,
$
    \proj_{T_{\cM(\overline{\z})}}(\v)
    \in
    \partial_{\cM(\overline{\z})}g(\overline{\z}).
$
It follows that
\[
\partial\left(g+\delta_{\cM(\overline{\z})}\right)(\overline{\z})
=
\partial_{\cM(\overline{\z})}g(\overline{\z})
+
T_{\cM(\overline{\z})}^{\perp}.
\]

Similarly, decomposing the Euclidean gradient into its tangential and normal
components gives
\[
    \nabla f(\overline{\z})
    =
    \proj_{T_{\cM(\overline{\z})}}\bigl(\nabla f(\overline{\z})\bigr)
    +
    \proj_{T_{\cM(\overline{\z})}^{\perp}}\bigl(\nabla f(\overline{\z})\bigr)
    =
    \nabla_{\cM(\overline{\z})}f(\overline{\z})
    +
    \proj_{T_{\cM(\overline{\z})}^{\perp}}\bigl(\nabla f(\overline{\z})\bigr).
\]
Substituting the last two identities into \eqref{eq:G_op}, and absorbing the
normal component of \(\nabla f(\overline{\z})\) into
\(T_{\cM(\overline{\z})}^{\perp}\), yields
\[
\partial G^t_{\gamma}(\y;\overline{\z})\mid_{\y=\overline{\z}}
=
\nabla_{\cM(\overline{\z})}f(\overline{\z})
+
\partial_{\cM(\overline{\z})}g(\overline{\z})
+
T_{\cM(\overline{\z})}^{\perp}.
\]
Now observe that
$
\nabla_{\cM(\overline{\z})}f(\overline{\z})
+
\partial_{\cM(\overline{\z})}g(\overline{\z})
\subseteq
T_{\cM(\overline{\z})},
$
whereas $T_{\cM(\overline{\z})}^{\perp}$ is orthogonal to
$T_{\cM(\overline{\z})}$. Hence,
$
\bz\in \partial G^t_{\gamma}(\overline{\z};\overline{\z})
$
holds if and only if the tangential component vanishes, i.e.,
\[
    \bz\in
    \nabla_{\cM(\overline{\z})}f(\overline{\z})
    +
    \partial_{\cM(\overline{\z})}g(\overline{\z}).
\]

	We complete the proof. 

\end{proof}

\subsection{Proof of Proposition \ref{pro:R}}
\label{sec:proof_conti}
\begin{proof}
We first prove the implication $(i)\Rightarrow(ii)$. Suppose that \(\overline R_\gamma^t(\overline{\z})=0\). 
By Proposition~\ref{pro:bridge}, this is equivalent to
$
    \overline T_\gamma^t(\overline{\z})=\overline{\z}.
$
Let $
    \{\z^k\}_{k\in\N}\subseteq \cX\cap\operatorname{int}(\dom h)
$ with $\z^k\to\overline{\z},$
be arbitrary, and set
$
    \y^k:=T_\gamma^t(\z^k).
$
Since $\z^k\in\operatorname{int}(\dom h)$, we have $
    T_\gamma^t(\z^k)=\overline T_\gamma^t(\z^k).
$
By Proposition~\ref{prop:extended-map-properties}~(ii), the sequence
$\{\y^k\}_{k\in\N}$ is bounded. Let $\overline{\v}$ be an arbitrary cluster
point of $\{\y^k\}_{k\in\N}$. Passing to a subsequence if necessary, we may
assume that $\y^k\to\overline{\v}$. By
Proposition~\ref{prop:extended-map-properties}~(iii), we obtain $    \overline{\v}_{\cB(\overline{\z})}
    =
    \overline{\z}_{\cB(\overline{\z})}$ and $\cB(\overline{\v})=\cB(\overline{\z})$. 
In particular, $\overline{\v}\in\cX\cap\cM(\overline{\z})$.

Using the optimality of $\y^k=T_\gamma^t(\z^k)$, with $\z^k$ as a comparison
point, gives
\begin{equation}
 \label{eq:prop_5.2_opt}   \gamma(\y^k;\z^k)+g(\y^k)+\frac1tD_h(\y^k,\z^k)
    \le
    \gamma(\z^k;\z^k)+g(\z^k).
\end{equation}
Since \(D_h\ge0\), we have
\[
    \gamma(\y^k;\z^k)+g(\y^k)
    +
    \frac1t
    \sum_{i\in\cI(\overline{\z})}
    D_\varphi(y_i^k,z_i^k)
    \le
    \gamma(\z^k;\z^k)+g(\z^k).
\]
Taking limits along the chosen subsequence, and using the continuity of
\(\gamma\), the local Lipschitz continuity of \(g\) on \(\cX\), and the
continuity of \(D_\varphi\) on the coordinates in \(\cI(\overline{\z})\), yields
\[
    \gamma(\overline{\v};\overline{\z})
    +
    g(\overline{\v})
    +
    \frac1t
    \sum_{i\in\cI(\overline{\z})}
    D_\varphi(\overline v_i,\overline z_i)
    \le
    f(\overline{\z})+g(\overline{\z}).
\]
Equivalently,
$
    G_\gamma^t(\overline{\v};\overline{\z})
    \le
    G_\gamma^t(\overline{\z};\overline{\z}).
$

Since $\overline T_\gamma^t(\overline{\z})=\overline{\z}$, the point $\overline{\z}$ is the unique minimizer of
$G_\gamma^t(\cdot;\overline{\z})$. Hence, we have 
$
    \overline{\v}=\overline{\z}.
$
Because every cluster point of $\{\y^k\}_{k\in\N}$ is $\overline{\z}$, we
obtain $
    T_\gamma^t(\z^k)=\y^k\to\overline{\z}.
$

Returning to the optimality inequality \eqref{eq:prop_5.2_opt},
we obtain
\[
    \frac1tD_h(\y^k,\z^k)
    \le
    \gamma(\z^k;\z^k)+g(\z^k)
    -
    \gamma(\y^k;\z^k)-g(\y^k).
\]
Since $\y^k\to\overline{\z}$, $\z^k\to\overline{\z}$, $\gamma$ is continuous,
and $g$ is locally Lipschitz continuous on $\cX$, we have
$
    \limsup_{k\to\infty}D_h(\y^k,\z^k)\le0.
$
Since $D_h\ge0$, it follows that
\[
    R_\gamma^t(\z^k)
    =
    D_h(T_\gamma^t(\z^k),\z^k)
    =
    D_h(\y^k,\z^k)
    \to0.
\]

Conversely, suppose that (ii) holds. We prove that $\overline R_\gamma^t(\overline{\z})=0$. 
By Proposition \ref{pro:bridge}, it suffices to prove $\overline{T}^t_{\gamma}(\overline{\z})=\overline{\z}$.
	We consider an interior point $\z^0\in \intm(\dom(h))\cap\cX$ and the sequence $\{\z^k\}_{k\in\N}$ defined by 
	\[\z^k:=\theta_k\z^0+(1-\theta_k)\overline{\z}, \quad\text{ with }\quad \theta_0=1,~ \theta_k\in(0,1],~\theta_k\to0.\]
Clearly, we have $\z^k\in \intm(\dom(h))\cap\cX$ and $\z^k\to\overline{\z}$, and thus the  limit condition ensures
	\[\lim_{k\to\infty}{R}^t_{\gamma}(\z^k)=\lim_{k\to\infty}D_h\left({T}^t_{\gamma}(\z^k),\z^k\right)=0. \]
It follows from Lemma \ref{lemma:bregman_limit} that the limit of $\{{T}^t_{\gamma}(\z^k)\}_{k\in \N}$ is $\overline{\z}$.

We now show that $\overline{\z}$ minimizes
$G_\gamma^t(\cdot;\overline{\z})$. Take any
$\y\in\cX\cap\cM(\overline{\z})
$
with $G_\gamma^t(\y;\overline{\z})<+\infty$; otherwise the desired inequality
is trivial. Define the recovery sequence
\[
    \widehat{\y}^k
    :=
    \theta_k\z^0+(1-\theta_k)\y.
\]
Then, we have $ \widehat{\y}^k\in\cX\cap\operatorname{int}(\dom h)$ and $\widehat{\y}^k\to\y$. 
Moreover, since $\y\in\cM(\overline{\z})$, we have $
    \y_{\cB(\overline{\z})}
    =
    \overline{\z}_{\cB(\overline{\z})}.
$
Hence
\[
    \widehat{\y}^k_{\cB(\overline{\z})}
    =
    \theta_k\z^0_{\cB(\overline{\z})}
    +
    (1-\theta_k)\overline{\z}_{\cB(\overline{\z})}
    =
    \z^k_{\cB(\overline{\z})}.
\]
Thus, $D_\varphi({\widehat{y}}_i^k,z_i^k)=0$ for all $i\in\cB(\overline{\z})$. 
Using the optimality of $\y^k=T_\gamma^t(\z^k)$, with $\widehat{\y}^k$ as a
comparison point, gives
\[
    \gamma(\y^k;\z^k)+g(\y^k)+\frac1tD_h(\y^k,\z^k)
    \le
    \gamma(\widehat{\y}^k;\z^k)+g(\widehat{\y}^k)
    +
    \frac1tD_h(\widehat{\y}^k,\z^k).
\]
Letting \(k\to\infty\), the left-hand side converges to
$
    f(\overline{\z})+g(\overline{\z})
    =
    G_\gamma^t(\overline{\z};\overline{\z}),
$ 
because $\y^k\to\overline{\z}$, $\z^k\to\overline{\z}$, and the condition (ii)
$D_h(\y^k,\z^k)\to0$. On the right-hand side, the Bregman terms over
$\cB(\overline{\z})$ vanish identically since
$\widehat y_i^k=z_i^k$ for all $i\in\cB(\overline{\z})$, while the terms over
$\cI(\overline{\z})$ converge by continuity.
Hence
\[
    \lim_{k\to\infty}D_h(\widehat{\y}^k,\z^k)
    =
    \sum_{i\in\cI(\overline{\z})}
    D_\varphi(y_i,\overline z_i).
\]
Consequently,
\[
    G_\gamma^t(\overline{\z};\overline{\z})
    \le
    \gamma(\y;\overline{\z})+g(\y)
    +
    \frac1t
    \sum_{i\in\cI(\overline{\z})}
    D_\varphi(y_i,\overline z_i)
    =
    G_\gamma^t(\y;\overline{\z}).
\]
Since $\y\in\cX\cap\cM(\overline{\z})$  is arbitrary, we conclude that
$\overline{\z}$ minimizes $G_\gamma^t(\cdot;\overline{\z})$. By the
uniqueness of the minimizer,
$
    \overline T_\gamma^t(\overline{\z})=\overline{\z}.
$
We complete the proof. 
\end{proof}

\section{Closing Remarks}\label{sec:conclusion}

This paper identifies a structural limitation of Bregman proximal-type methods under singular Bregman geometries. We show that spurious stationary points can arise systematically, that standard Bregman residuals may certify only restricted stationarity near the boundary, and that the resulting small Bregman displacements can lead to arbitrarily long finite-time trapping. These findings call for a more careful interpretation of residual-based convergence guarantees and motivate alternative certificates or algorithmic safeguards.



\section*{Acknowledgments}
Jiajin Li thanks Prof. Heinz H. Bauschke, Prof. Joseph Paat and  Prof. Yinyu Ye  for their helpful discussions. Jiajin Li was supported by a Natural Sciences and Engineering Research Council of Canada Discovery Grant RGPIN-2025-05817.
\bibliography{ref}
\bibliographystyle{abbrvnat}

\appendix
\section{Verification of Assumption \ref{assum:gamma} 
	(iv)}\label{appen:verify}
In this section, we verify that Assumption~\ref{assum:gamma}~(iv) is satisfied by several standard surrogate models commonly used in the literature.
\begin{proposition}\label{prop:verify1}
	Suppose that Assumption~\ref{assum:h} holds. Then the following statements hold:
	\begin{enumerate}[label={\rm (\roman*)}]
		\item If the surrogate model takes the form $\gamma(\y;\x)=f(\x)+\nabla f(\x)^{\top}(\y-\x)$ and      $h+tg$ is supercoercive for all $t\in(0,\overline{t}]$ for some $\overline{t}>0$, then Assumption \ref{assum:gamma} (iv) is satisfied.
		\item If the surrogate model takes the form $\gamma(\y;\x)=f(\y)$, $\inf_{\x} F(\x)>-\infty$, and the set
    $
        \cD(\z,\alpha):=\{\y\in\dom(h):D_h(\y,\z)\leq \alpha\}
    $
    is bounded for every $\z\in\intm(\dom(h))$ and every $\alpha\in\R$, then Assumption \ref{assum:gamma} (iv) is satisfied.
		\item If the surrogate model  takes the form 
		$\gamma(\y;\x)=f(\x)+\nabla f(\x)^{\top}(\y-\x)+\frac12(\y-\x)^{\top}\nabla^2 f(\x)(\y-\x)$, where $f$ is twice continuously differentiable and convex, and if
$h+tg$ is supercoercive for all $t\in(0,\overline t]$ for some
$\overline t>0$, then Assumption~\ref{assum:gamma}~(iv) is satisfied.
	\end{enumerate}
\end{proposition}
\begin{remark}
	The assumptions made in Proposition~\ref{prop:verify1} are consistent with standard practices in the literature:
	(i) For the Bregman proximal gradient method with surrogate model $\gamma(\y;\x) = f(\x) + \nabla f(\x)^{\top}(\y - \x)$, the supercoercivity of $h + tg$ is a common assumption; see, e.g., \citep[Lemma 2]{bauschke2017descent} and \citep[Assumption B]{bolte2018first}. 
	(ii) The lower boundedness of $F$ and level boundedness of $D_h(\cdot,\x)$ are standard assumptions in the analysis of Bregman proximal point methods; see, e.g., \citep{chen1993convergence,yang2022bregman}.
	(iii) The convexity condition on $f$ in Proposition~\ref{prop:verify1}~(iii) is also assumed in \citet{doikov2023gradient},  and the supercoercive property of $h+tg$ is satisfied for all examples therein.
\end{remark}

\begin{proof}
(i): Fix $t\in(0,\overline{t}]$.
To verify Assumption~\ref{assum:gamma} (iv), we prove the stronger statement:
	\[\lim\limits_{k\to\infty}\frac{\gamma(\y^k;\z^k)+g(\y^k)+\frac1t D_h(\y^k,\z^k)}{\|\y_k\|}=+\infty.
	\]
Since $\z^k\to\overline{\z}$ and $\nabla f$ is continuous, the sequence
$\{\nabla f(\z^k)\}_{k\in\N}$ is bounded. Moreover, the sequences
$\{f(\z^k)\}_{k\in\N}$ and
$\{\nabla f(\z^k)^\top \z^k\}_{k\in\N}$ are bounded. Hence, there exists
$C>0$ such that
\[
\begin{aligned}
    \gamma(\y^k;\z^k)
    &= f(\z^k)+\nabla f(\z^k)^\top(\y^k-\z^k) \ge -C(1+\|\y^k\|).
\end{aligned}
\]
Therefore, it suffices to show that
\[
    \lim_{k\to\infty}
    \frac{t g(\y^k)+D_h(\y^k,\z^k)}{\|\y^k\|}
    =+\infty . 
\]
Since $\|\y^k\|\to+\infty$ and $\y^k\in\intm(\dom h)$, the interval
$\dom(\varphi)$ must be unbounded. We only prove the case where
$\cl(\dom(\varphi))=[a,+\infty)$ with finite $a$; the other unbounded cases
are handled analogously.


{
Since $\z^k\to\overline{\z}$, there exists 
$\eta\in(a,+\infty)$ such that
\[
    z_i^k\le \eta,\qquad \forall i\in[n],\ k\in\N .
\]
Define the large-coordinate index set
$
    \cJ(\y^k):=\{j\in[n]: y_j^k\ge \eta\}.
$
Then, we have

\begin{equation*}
\label{eq:super1}
	\begin{aligned}
	\lim\limits_{k\to\infty}\frac{t g(\y^k) + D_h(\y^k,\z^k)}{\|\y^k\|}&\geq \lim\limits_{k\to\infty}\frac{tg(\y^k) +\sum_{j\in\cJ(\y^k)} D_{\varphi}(y^k_j,z^k_j)}{\|\y^k\|}\\
	&\geq \lim\limits_{k\to\infty}\frac{tg(\y^k) +\sum_{j\in\cJ(\y^k)} D_{\varphi}(y^k_j,\eta)}{\|\y^k\|}\\
	&= \lim\limits_{k\to\infty}\frac{tg(\y^k) +\sum_{j\in\cJ(\y^k)} \varphi(y^k_j)-\varphi(\eta)-\varphi^{\prime}(\eta)(y^k_j-\eta) }{\|\y^k\|}\\
	&\geq \lim\limits_{k\to\infty}\frac{tg(\y^k) +\sum_{j\in\cJ(\y^k)} \varphi(y^k_j)}{\|\y^k\|}-n|\varphi^{\prime}(\eta)|,
	\end{aligned}
\end{equation*}
where the second inequality is due to the inequality \eqref{eq:three-point-identity}  and $y^k_j\geq\eta\geq z^k_j$ for $j\in\cJ(\y^k)$, and the last inequality follows from 
$
    \sum_{j\in\cJ(\y^k)} |y_j^k|\le n\|\y^k\|. 
$

We now show that the limit on the right-hand side is $+\infty$.
Since
\(|\cJ(\y^k)|\le n\), we have
\[
    \frac{tg(\y^0)+\sum_{j\in\cJ(\y^k)}\varphi(y_j^0)}{\|\y^k\|}
    \to 0.
\]
Therefore, by the convexity of \(\varphi\) and \(g\),
\begin{equation}\label{eq:super2}
\begin{aligned}
\lim_{k\to\infty}
\frac{tg(\y^k)+\sum_{j\in\cJ(\y^k)}\varphi(y_j^k)}{\|\y^k\|}
&=
\lim_{k\to\infty}
\frac{
tg(\y^k)+tg(\y^0)
+\sum_{j\in\cJ(\y^k)}
\left(\varphi(y_j^k)+\varphi(y_j^0)\right)
}{\|\y^k\|} \\
&\ge
\lim_{k\to\infty}
2\cdot
\frac{
tg\left(\frac{\y^k+\y^0}{2}\right)
+\sum_{j\in\cJ(\y^k)}
\varphi\left(\frac{y_j^k+y_j^0}{2}\right)
}{\|\y^k\|}.
\end{aligned}
\end{equation}
Set
$
    \beta_1:=\min_{i\in[n]}y_i^0$ and $
    \beta_2:=\max_{i\in[n]}y_i^0$.
 Note that for $j\notin\cJ(\y^k)$, we have $y_j^k<\eta$. Hence
\[
    \frac{y_j^k+y_j^0}{2}
    \in
    \left[\frac{a+\beta_1}{2},\frac{\eta+\beta_2}{2}\right]
    \subset \intm(\dom\varphi),
    \qquad \forall j\notin\cJ(\y^k),\ k\in\N .
\]
By the continuity of $\varphi$, there exists $M>0$ such that
$
|
    \varphi\left(\frac{y_j^k+y_j^0}{2}\right)
  |
    \le M, \forall j\notin\cJ(\y^k),\ k\in\N .
$
It follows that
\[
    \lim_{k\to\infty}
    \frac{
        \sum_{j\notin\cJ(\y^k)}
        \varphi\left(\frac{y_j^k+y_j^0}{2}\right)
    }{\|\y^k\|}
    =0 .
\]
Combining this with \eqref{eq:super2}, we obtain
\[
\begin{aligned}
\lim_{k\to\infty}
\frac{
    t g(\y^k)+\sum_{j\in\cJ(\y^k)}\varphi(y_j^k)
}{\|\y^k\|}
&\ge
\lim_{k\to\infty}
2\cdot
\frac{
    t g\left(\frac{\y^k+\y^0}{2}\right)
    +
    h\left(\frac{\y^k+\y^0}{2}\right)
}{\|\y^k\|}.
\end{aligned}
\]
Finally, since
$
    \|\frac{\y^k+\y^0}{2}\|
    \to+\infty$ 
and $
    \frac{2\left\|\frac{\y^k+\y^0}{2}\right\|}{\|\y^k\|}\to 1,
$
the supercoercivity of $tg+h$ yields
\[
\lim_{k\to\infty}
\frac{
    t g(\y^k)+\sum_{j\in\cJ(\y^k)}\varphi(y_j^k)
}{\|\y^k\|} \ge \lim_{k\to\infty}
2\cdot
\frac{
    t g\left(\frac{\y^k+\y^0}{2}\right)
    +
    h\left(\frac{\y^k+\y^0}{2}\right)
}{\|\y^k\|}
=+\infty .
\]
Putting everything together proves (i).

}

(ii) Fix $t\in(0,\bar t]$.
Since $\gamma(\y^k;\z^k)=f(\y^k)$, we have
	\[\begin{aligned}
		\lim\limits_{k\to\infty} \gamma(\y^k;\z^k)+g(\y^k)+\frac1t D_h(\y^k,\z^k)&=\lim\limits_{k\to\infty} F(\y^k)+\frac1t D_h(\y^k,\z^k)\\
		&\geq \inf_{\x}F(\x)+\lim\limits_{k\to\infty}\frac1t D_h(\y^k,\z^k).
	\end{aligned} \]
	Due to $\inf_{\x}F(\x)>-\infty$, to verify Assumption~\ref{assum:gamma} (iv), it suffices to show that
	\[\lim\limits_{k\to\infty} D_h(\y^k,\z^k)=+\infty. \]
As before, we only prove the case where
$\cl(\dom(\varphi))=[a,+\infty)$ with finite $a$; the other unbounded cases
are analogous.

Since $\z^k\to\overline{\z}$, we can choose
$\eta\in(a,+\infty)$ such that
$
    z_i^k\le \eta, \forall i\in[n],\ k\in\N .
$
Define $\widehat{\y}^k$ coordinatewise by
$
    \widehat y_i^k:=\max\{y_i^k,\eta\},\forall  i\in[n].
$
Then $\widehat{\y}^k\in\intm(\dom h)$ and
$\|\widehat{\y}^k\|\to+\infty$. For each $i\in[n]$, we compare the one-dimensional Bregman divergences.
If $y_i^k\ge \eta\ge z_i^k$, then \eqref{eq:three-point-identity} gives
\[
    D_\varphi(y_i^k,z_i^k)
    \ge
    D_\varphi(y_i^k,\eta)
    =
    D_\varphi(\widehat y_i^k,\eta).
\]
If $y_i^k\le \eta$, then $\widehat y_i^k=\eta$, and hence
$
    D_\varphi(y_i^k,z_i^k)
    \ge 0
    =
    D_\varphi(\widehat y_i^k,\eta).
$
Therefore,
\[
    D_\varphi(y_i^k,z_i^k)
    \ge
    D_\varphi(\widehat y_i^k,\eta),
    \qquad \forall i\in[n],\ k\in\N.
\]
It follows that
\[\lim\limits_{k\to\infty} D_h(\y^k,\z^k)\geq \lim\limits_{k\to\infty} D_h(\widehat{\y}^k,\eta\1)=+\infty.  \]
Here, the equality is due to $\|\widehat{\y}^k\|\to+\infty$, $\eta\1\in\intm(\dom(h))$, and the level boundedness condition, which ensures the boundedness of $\cD(\eta\1,\alpha)$ for all $\alpha\in\R$. We complete the proof of (ii). 

(iii) Fix $t\in(0,\overline t]$. Since $f$ is twice continuously
differentiable and convex,
$
    (\y-\x)^\top\nabla^2 f(\x)(\y-\x)\ge0, \forall \x,\y\in\R^n .
$
Thus the second-order surrogate satisfies
$
    \gamma(\y^k;\z^k)
    \ge
    f(\z^k)+\nabla f(\z^k)^\top(\y^k-\z^k).
$
Therefore,
\[
\gamma(\y^k;\z^k)+g(\y^k)+\frac1tD_h(\y^k,\z^k)
\ge
f(\z^k)+\nabla f(\z^k)^\top(\y^k-\z^k)
+g(\y^k)+\frac1tD_h(\y^k,\z^k),
\]
and the right-hand side tends to \(+\infty\) by part (i). This completes the
proof.
	
\end{proof}

\section{Missing Proofs for Example \ref{example:nonconvex}} 
\label{appen:nonconvex_spurious}
The omitted technical details for Example~\ref{example:nonconvex} are provided here. We first verify that the function defined in~\eqref{eq:concrete_phi} satisfies the assumptions stated in Example~\ref{example:nonconvex}.
\begin{fact}\label{fact:phi}
	Suppose that $\alpha \in (0,0.1]$ and consider the Burg entropy kernel $\varphi(x) = -\log(x)$. The function $f_{\alpha}:\R^2 \rightarrow \R$ defined in Example~\ref{example:nonconvex}, with $\phi_\alpha$ given by~\eqref{eq:concrete_phi}, satisfies the following properties:
	\begin{enumerate}[label={{\rm (\roman*)}}, itemsep=0.5pt]
		\item $\phi_{\alpha}$ is continuously differentiable; 
        \item $\phi^{\prime}_{\alpha}(x)\geq1$ for $x\in[0,1]$;
		\item $\phi_{\alpha}(x)=2(x-\alpha)$ for $x\geq \alpha$, $\phi_{\alpha}(x)\leq0$ for $x\leq \alpha$, and $\phi_{\alpha}(0)\leq-1$;
		\item  $f_{\alpha}+h$ and $-f_{\alpha}+h$ are convex on $(0,1]\times(0,1]$.
	\end{enumerate}
\end{fact}
\begin{proof}
	(i): 
	Establishing differentiability over $\R$ reduces to verifying both continuity and differentiability of $\phi_\alpha$ at the junction points where the functional form changes. These points are:
	\[
	x = \alpha,\quad x = \alpha \exp\left(-\frac{1}{\alpha}\right),\quad \text{and} \quad x = \frac{1}{2}\alpha \exp\left(-\frac{1}{\alpha}\right).
	\]
	As a representative case, we verify continuity and differentiability at 
	$x = \alpha$:
	\begin{itemize}[itemsep=0.5pt]
		\item For $x \ge \alpha$, we have $\phi_\alpha(x) = 2x - 2\alpha$, which gives $\phi_\alpha(\alpha^+) = 0$ and $\phi_\alpha'(\alpha^+) = 2$.
		\item For $x \in \left[\alpha \exp\left(-\tfrac{1}{\alpha}\right), \alpha\right]$, we have $\phi_\alpha(x) = x + \alpha \log\left(\tfrac{x}{\alpha}\right) - \alpha$. Hence, we have 
		\[
		\phi_\alpha(\alpha^-) = \alpha + \alpha \log(1) - \alpha = 0, \quad \text{and } \quad 
		\phi_\alpha'(\alpha^-) = 1 + \frac{\alpha}{x} \Big|_{x = \alpha} = 2.
		\]
	\end{itemize}
	Thus, $\phi_\alpha$ is continuously differentiable at $x = \alpha$. Similar computations at the other two junction points confirm that $\phi_\alpha$ and its derivative are continuous at those points as well.

    (ii): We verify that $\phi_\alpha'(x) \ge 1$ for all $x \in [0,1]$ by examining the derivative on each interval specified in~\eqref{eq:concrete_phi}.
	\begin{itemize}[itemsep=0.5pt]
		\item For $x \in [\alpha, 1]$, we have $\phi_\alpha'(x) = 2$.
		\item For $x \in \left[\alpha\exp\left(-\tfrac{1}{\alpha}\right), \alpha\right]$, we have $\phi_\alpha'(x) = 1 + \frac{\alpha}{x} \ge 2$. 
		\item  For $x \in \left[\tfrac{1}{2} \alpha\exp\left(-\tfrac{1}{\alpha}\right), \alpha\exp\left(-\tfrac{1}{\alpha}\right)\right]$, we have $\phi_\alpha'(x) = 1-\frac{\alpha}{x} + 2 \exp\left(\tfrac{1}{\alpha}\right) \ge 1$.
		\item  For $x \in \left[0, \tfrac{1}{2} \alpha\exp\left(-\tfrac{1}{\alpha}\right)\right]$,  we have $\phi_\alpha'(x) = 1$.
	\end{itemize}
	Therefore, $\phi_\alpha'(x) \ge 1$ for all $x \in [0,1]$.
	
	(iii): This follows directly from the definition of $\phi_\alpha(x)$ in~\eqref{eq:concrete_phi}:
	\begin{itemize}[itemsep=0.5pt]
		\item For $x \ge \alpha$, $\phi_\alpha(x) = 2(x - \alpha)$ by construction.
		\item  By $\phi_\alpha'(x) \ge 1$ for $x \in [0,1]$, $\phi_\alpha$ monotonically increases on $[0,1]$. Hence, $\phi_\alpha(x)\leq \phi_\alpha(\alpha)=0$ for $x \le \alpha$. 
		\item In particular, $\phi_\alpha(0) = \alpha(\log 2 - 2) - 1 \le -1$ for all $\alpha \in (0, 0.1]$.
	\end{itemize}
	
	(iv): We begin by verifying that, on each piece of 
	$f_{\alpha}$, the Hessians of $h + f_{\alpha}$ and $h - f_{\alpha}$ are positive semidefinite.
	To this end, we examine each case of 
	$\phi_\alpha(x_1)$ piece by piece:
	
	\begin{itemize}[itemsep=0.5pt]
		\item For $x_1 \in [\alpha, 1]$, we have 
		\[\frac{1}{x_1^2}-| \nabla^2_{11} f_{\alpha}(\x)|=\frac{1}{x_1^2}\geq |\nabla^2_{12} f_{\alpha}(\x)| =1, \text{ and  } \frac{1}{x_2^2}-| \nabla^2_{22} f_{\alpha}(\x)|=\frac{1}{x_2^2}\geq |\nabla^2_{21} f_{\alpha}(\x)| =1.\]
		This implies that $\nabla^2 h+\nabla^2f_\alpha$ and  $\nabla^2 h-\nabla^2f_\alpha$ are diagonally dominant with strictly positive diagonal entries, and are therefore positive semidefinite.
		\item For $x_1 \in \left[\alpha\exp\left(-\tfrac{1}{\alpha}\right), \alpha\right]$, we have 
		\[\frac{1}{x_1^2} -| \nabla_{11}^2f_{\alpha}(\x)|=\frac{1}{x_1^2}-\frac{\alpha}{2}\cdot \frac{1}{x_1^2}(x_2+0.05)\geq \frac{0.9}{x_1^2}>0.\] 
		
		Using $\frac{1}{x_1^2} -| \nabla_{11}^2f_{\alpha}(\x)|\geq\frac{0.9}{x_1^2}$ and $\frac{1}{x_2^2} -|\nabla_{22}^2f_{\alpha}(\x)|=\frac1{x_2^2}\geq1$, we have
		\[\left(\frac{1}{x_1^2} -| \nabla_{11}^2f_{\alpha}(\x)|\right)\left(\frac{1}{x_2^2} -|\nabla_{22}^2f_{\alpha}(\x)|\right)\geq \frac{0.9}{x_1^2}\geq |\nabla^2_{12} f_{\alpha}(\x)|\cdot|\nabla^2_{21} f_{\alpha}(\x)| =\frac14\left(1+\frac{\alpha}{x_1}\right)^2,\]
         where the last inequality is due to
        \[\sqrt{\frac{0.9}{x_1^2}}\geq\frac{0.8}{x_1}=\frac12\left(\frac{1}{x_1}+\frac{0.6}{x_1}\right)\geq\frac{1}{2}\left(1+\frac{\alpha}{x_1}\right).\] 
        Hence, all leading principal minors of  $\nabla^2 h+\nabla^2f_\alpha$ and  $\nabla^2 h-\nabla^2f_\alpha$ are positive, which implies that both matrices are positive definite.
		
		\item  For $x_1 \in  \left[\tfrac{1}{2} \alpha\exp\left(-\tfrac{1}{\alpha}\right), \alpha\exp\left(-\tfrac{1}{\alpha}\right)\right]$, we have 
		\[\frac{1}{x_1^2} -|\nabla^2_{11} f_{\alpha}(\x)|=\frac{1}{x_1^2}-\frac{\alpha}{2}\cdot \frac{1}{x^2_1}(x_2+0.05)\geq \frac{0.9}{x_1^2}>  0.\]
		Moreover, using the estimates
		\[\frac{1}{x_1^2} -| \nabla_{11}^2f_{\alpha}(\x)|\geq\frac{0.9}{x_1^2};\qquad \frac{1}{x_2^2} -|\nabla_{22}^2f_{\alpha}(\x)|=\frac1{x_2^2}\geq1;\]
		\[ |\nabla^2_{12} f_{\alpha}(\x)|=|\nabla^2_{21} f_{\alpha}(\x)|= \frac12\left(1-\frac{\alpha}{x_1}+2\exp\left(\frac{1}{\alpha}\right)\right)\leq\frac12+\exp\left(\frac1{\alpha}\right)\leq2\exp\left(\frac1{\alpha}\right),\]
		we have
		\[\left(\frac{1}{x_1^2} -| \nabla_{11}^2f_{\alpha}(\x)|\right)\left(\frac{1}{x_2^2} -|\nabla_{22}^2f_{\alpha}(\x)|\right)\geq \frac{0.9}{x_1^2}>\left(2\exp\left(\frac1{\alpha}\right)\right)^2\geq|\nabla^2_{12} f_{\alpha}(\x)|\cdot|\nabla^2_{21} f_{\alpha}(\x)| .\] 
		Here, the second inequality follows from 
		$x_1 \in \left[\tfrac{1}{2} \alpha\exp\left(-\tfrac{1}{\alpha}\right), \alpha\exp\left(-\tfrac{1}{\alpha}\right)\right]$ and $\alpha\leq0.1$, which together imply that
		\[
		\frac{0.9}{x_1^2} \ge \frac{0.9 \exp\left(\frac{2}{\alpha}\right)}{\alpha^2} \ge 90\exp\left(\frac{2}{\alpha}\right) > \left(2\exp\left(\frac{1}{\alpha}\right)\right)^2. 
		\]
		
		This confirms that both $\nabla^2 h+\nabla^2f_\alpha$ and  $\nabla^2 h-\nabla^2f_\alpha$ have positive leading principal minors, and thus are positive definite. 
		
		\item  For $x_1 \in \left[0, \tfrac{1}{2} \alpha\exp\left(-\tfrac{1}{\alpha}\right)\right]$, the same reasoning as in the first case applies. 
	\end{itemize}
	
	Consequently, both $h + f_{\alpha}$ and $h - f_{\alpha}$ are piecewise convex. Since both of them are also continuously differentiable across pieces, the conditions of \citet[Theorem 5.5]{bauschke2016convexity} are satisfied; see also the verification argument from \citet[Example 6.1]{bauschke2016convexity}. We thus conclude that $h + f_{\alpha}$ and $h - f_{\alpha}$ are convex on $(0,1]\times (0,1]$.
	
This completes the proof.
\end{proof}

Finally, we present the detailed proof of Observation~\ref{observation:encounter2}.
\begin{proof}[Proof of Observation \ref{observation:encounter2}]
	We set the parameter $\alpha$ and the index $k_1$ as follows:
	\[\alpha=\frac{2}{\ceil{40\exp\left(\frac1{\epsilon}\right)}+K}, \qquad k_1= \ceil{40\exp\left(\frac{1}{\epsilon}\right)}.\]
	Let $k_2$ be the minimal index such that $x_1^k < 2\alpha$, i.e., $k_2=\min\{k\in\N:x_1^k< 2\alpha\}$. Our first goal is to prove that  $k_1+K< k_2$.  
	
	A direct computation shows that $\nabla_{x_1} f_\alpha(\x^k)>0$. Combined with the KKT conditions \eqref{eq:kkt_nonconvex}, this implies that $x_1^{k+1}<x_1^k\leq 1, \lambda_1^k=0$, and 
	\begin{equation}\label{eq:encounter_update1}
		\frac{1}{x^{k+1}_1}=\frac{1}{x^k_1} + \frac14 \phi_{\alpha}^{\prime}(x_1^k)\cdot(x_2^k+0.05).
	\end{equation}
	For the second coordinate $x_2$, we recall the update rule from \eqref{eq:x2_update} and $t=0.5$: 
		\begin{equation}\label{eq:encounter_update2}
		\frac{1}{x_2^{k+1}} = \max\left(\frac{1}{x_2^k} + \frac{1}{4} \phi_\alpha(x_1^k),\ 1\right).
	\end{equation}
	By the definition of $k_2$, we have $x^k_1\geq2\alpha$ for $k< k_2$. It then follows from the definition of $\phi_\alpha$ that $\phi_\alpha^\prime(x_1^k)=2$ and $\phi_\alpha(x_1^k)\ge 0$ for all $k<k_2$. Combining these with the update rule \eqref{eq:encounter_update1} and \eqref{eq:x2_update}, we obtain $x_2^k\leq x_2^0 = 0.1$ and 
	\begin{equation}\label{eq:k2}
		\frac{1}{x_1^k} 
		< \frac{1}{x_1^{k-1}} + 0.1 
		< \cdots 
		< \frac{1}{x_1^0} + 0.1 \cdot k 
		= 1 + \frac{k}{10}, 
		\qquad \forall\, k \leq k_2.
	\end{equation}
	Then, armed with the fact that $x_1^{k_2}<2\alpha$ and \eqref{eq:k2}, we get 
	\[
	\frac{1}{x_1^{k_2}} >\frac{1}{2\alpha} \quad \text {  and   } \quad \frac{1}{x_1^{k_2}} <1+\frac{k_2}{10}. 
	\]
	Combining the two gives $k_2>10\cdot (\tfrac{1}{2\alpha}-1)$. Since $\alpha\leq 0.1$, it follows that 
	\[
	\frac{1}{2\alpha} - 1 
	= \frac{1}{2\alpha} \cdot (1 - 2\alpha) 
	\geq \frac{1}{2\alpha} \cdot (1 - 0.2) 
	= \frac{0.8}{2\alpha} 
	= \frac{4}{10\alpha} \quad \text{and } \quad  k_2 >\frac{4}{\alpha}. 
	\]
	By the choice of $\alpha$, we have
$
    k_1+K=\frac{2}{\alpha}.
$
Since $k_2>4/\alpha$, it follows that
$
    k_1+K=\frac{2}{\alpha}<\frac{4}{\alpha}<k_2.
$
	
	Next, we are ready to continue to control $\|\x^k-\tilde{\x}^\star\|$ for $k\in[k_1,k_1+K]$. Combining the update rule \eqref{eq:encounter_update2} with the fact that $\phi_\alpha(x) = 2(x-\alpha)$ for $x\ge \alpha$, and using the bound $x_1^k\ge 2\alpha $ for all $k<k_2$, we obtain the following estimate for $x_2^k$: 
	\[
	\frac{1}{x_2^k}
	\geq \frac{1}{x_2^{k-1}} + \frac{1}{2}(x_1^{k-1} - \alpha)
	\geq \frac{1}{x_2^{k-1}} + \frac{1}{4} x_1^{k-1}
	\geq \cdots
	\geq \frac{1}{x_2^0} + \frac{1}{4} \sum_{i=0}^{k-1} x_1^i,
	\qquad \forall~ k \leq k_2.
	\]
	Using $x_2^0 = 0.1$ and the bound from \eqref{eq:k2}, i.e., $x_1^i\ge \tfrac{10}{i+10}$, we further obtain 
	\[  \frac{1}{x^k_2} \geq 10+\frac14\sum_{i=0}^{k-1}\frac{10}{i+10}\geq10+\frac52\sum_{i=0}^{k-1}\int_i^{i+1}\frac{1}{x+10} {\rm d}x= 10+\frac52\log\left(1+\frac{k}{10}\right), \qquad\ \forall\ k\leq k_2.
	\]
  Note that $k_1+K<k_2$, and $\log(1+\tfrac{k}{10})\ge \frac{1}{\epsilon}$ for $k\geq k_1= \ceil{40\exp\left(\frac{1}{\epsilon}\right)}$.
	The bound above implies 
	\[
	x_2^k \leq  \left( 10 + \frac{5}{2} \cdot \frac{1}{\epsilon} \right)^{-1}
	= \frac{1}{10 + \frac{5}{2} \cdot \frac{1}{\epsilon}}
	= \frac{2\epsilon}{20\epsilon + 5}
	\leq \frac{2}{5} \epsilon< \frac{1}{2}\epsilon,
	\qquad \forall\, k \in [k_1, k_1 + K].
	\]
	
	On the other hand, using the update rule \eqref{eq:encounter_update1}, along with the facts that  $x_2^k>0$ and $\phi_\alpha^\prime(x^k_1)=2$ for all $k<k_2$, we obtain
	\[ \frac{1}{x^k_1}\geq \frac{1}{x^{k-1}_1}+0.025\geq\cdots\geq \frac{1}{x^0_1}+\frac{k}{40}\geq 1+\frac{k_1}{40}>\exp\left(\frac{1}{\epsilon}\right)>\frac{2}{\epsilon},\qquad\forall\ k\in[k_1,k_1+K], \]
	where the last inequality holds because $\epsilon\in(0,1]$ and the function $x\mapsto\exp(x)-2x$ is positive on $[1,+\infty)$. It follows that
	\[x_1^k\leq\frac{\epsilon}{2}, \qquad\quad\forall\  k\in [k_1,k_1+K].\] 
	Combining this with $\tilde{\x}^\star=\bz$ and the earlier bound $x_2^k \leq \tfrac{\epsilon}{2}$ for $ k\in [k_1,k_1+K]$, we conclude that 
	\[\|\x^k-\tilde{\x}^\star\|\leq \epsilon,  \qquad\quad\forall\  k\in [k_1,k_1+K].\] 
	This completes the proof.
\end{proof}

\section{Proof of Proposition \ref{prop:extended-map-properties}}
\label{app:extended-update}
We first prove a coerciveness property of $G_\gamma^t$, which will be used to
establish well-posedness and boundedness.

\begin{lemma}[Coerciveness of $G$]\label{lemma:G}
	Suppose that the sequence $\{\z^k\}_{k\in\N}\subseteq\cX$  converges to $\overline{\z}\in\cX$ and the sequence $\{\y^k\}_{k\in\N}\subseteq\cX$ satisfies $\|\y^k\|\to+\infty$. Then, we have
	\[\lim_{k\to\infty}G^t_{\gamma}(\y^k;\z^k)=+\infty.\]
\end{lemma}
\begin{proof}

WLOG, we may assume that $G_\gamma^t(\y^k;\z^k)<+\infty$ for all $k\in\N$, since
otherwise the assertion is trivial along the indices where
$G_\gamma^t(\y^k;\z^k)=+\infty$. Then the indicator term in
$G_\gamma^t(\y^k;\z^k)$ enforces
$
    \y^k_{\cB(\z^k)}\equiv \z^k_{\cB(\z^k)}, \forall k\in\N .
$
	
Let $\z^{\intm}  \in\intm(\dom(h))\cap\cX$ be any interior point. For each $k\in\N$, we define 
	\[
	\widehat{\z}^k=(1-\theta_k)\z^k+\theta_k\z^{\intm}, \quad \widehat{\y}^k=(1-\theta_k)\y^k+\theta_k\z^{\intm},
	\]
	where  $\theta_k\in(0,1)$ with $\theta_k\to 0$. Then, it is easy to verify that $\widehat{\z}^k,\widehat{\y}^k\in\intm(\dom(h))\cap\cX$, $\widehat{\z}^k\to\overline{\z}$, $\|\widehat{\y}^k\|\to+\infty$, and $\widehat{\y}^k_{\cB(\z^k)}= \widehat{\z}^k_{\cB(\z^k)}$ for all $k \in \N$. It implies that 
	\[	\gamma(\widehat{\y}^k;\widehat{\z}^k)+g(\widehat{\y}^k)+\frac1t D_h(\widehat{\y}^k,\widehat{\z}^k)= \gamma(\widehat{\y}^k;\widehat{\z}^k)+g(\widehat{\y}^k)+\frac1t \sum_{i\in \cI(\z^k)}D_{\varphi}(\widehat{y}^k_i,\widehat{z}^k_i).
	\]
For each fixed $k$, by the continuity of $\gamma,g,D_\varphi$, we may choose
$\theta_k>0$ sufficiently small, with $\theta_k\downarrow 0$, such that
\begin{equation*} 
\left|\gamma(\widehat{\y}^k;\widehat{\z}^k)+g(\widehat{\y}^k)+\frac1t D_h(\widehat{\y}^k,\widehat{\z}^k)-G^t_{\gamma}(\y^k;\z^k)\right|\leq 1.
\end{equation*}
Now, by Assumption \ref{assum:gamma} (iv), we know $\gamma(\widehat{\y}^k;\widehat{\z}^k)+g(\widehat{\y}^k)+\frac1t D_h(\widehat{\y}^k,\widehat{\z}^k) \to +\infty$. Consequently, we can conclude that $G^t_{\gamma}(\y^k;\z^k)\to +\infty$. 
This completes the proof.
\end{proof}

With Lemma \ref{lemma:G} in place, we give  the full proof of Proposition~\ref{prop:extended-map-properties}.
\begin{proof}
	(i) In view of Assumption~\ref{assum:gamma} (iii), the function  $G^t_{\gamma}(\ \cdot\ ;\x)$ is strictly convex for any $\x\in\cX$. To show the well-posedness of the extended Bregman update mapping, it suffices to verify that $G^t_{\gamma}(\ \cdot\ ;\x)$ is level bounded. Lemma~\ref{lemma:G} establishes the coerciveness of $G^t_{\gamma}(\ \cdot\ ;\x)$  for any fixed $\x\in\cX$.  Moreover, $G_\gamma^t(\cdot;\x)$ is lower
semicontinuous, since $\gamma(\cdot;\x)$ is continuous, $g$ is lower
semicontinuous, the Bregman terms are lower semicontinuous, and
$\cM(\x)$ is closed.
 Hence, the level set $\{\y\in\cX: G^t_{\gamma}(\y ;\x)\leq c \}$ is compact for every $c\in\R$, and the minimizer exists and is unique. 

(ii) Let $\{\z^k\}_{k\in\N}\subset \cX$ be a bounded sequence. 
Suppose, to the contrary, that 
$\{\overline{T}^{t}_{\gamma}(\z^k)\}_{k\in\N}$ is unbounded. 
Passing to a subsequence if necessary, we may assume that $
    \|\overline{T}^{t}_{\gamma}(\z^k)\|\to+\infty.$ 
Since $\{\z^k\}_{k\in\N}$ is bounded and $\cX$ is closed, by passing to a
further subsequence if necessary, we may assume that
$
    \z^k\to\overline{\z}\in\cX .
$ 
Since
$
    G_\gamma^t(\z^k;\z^k)
    =
    \gamma(\z^k;\z^k)+g(\z^k),
$
Assumption~\ref{assum:gamma}~(i)--(ii) and the continuity of $g$ on $\cX$ yield
\[
	\lim_{k\to+\infty} G^t_{\gamma}(\z^k;\z^k)= \gamma(\overline{\z};\overline{\z})+g(\overline{\z})= f(\overline{\z}) + g(\overline{\z}) <+\infty. 
	\]
By the optimality of $\overline{T}^{t}_{\gamma}(\z^k)$, we have $
    G_\gamma^t(\overline{T}^{t}_{\gamma}(\z^k);\z^k)
    \le
    G_\gamma^t(\z^k;\z^k).$ 
Hence
\[
	\limsup_{k\to+\infty}G^t_{\gamma}(\overline{T}^{t}_{\gamma}(\z^k);\z^k) \leq   \lim_{k\to+\infty} G^t_{\gamma}(\z^k;\z^k)<+\infty. 
	\]
This contradicts Lemma~\ref{lemma:G}, applied with
$\y^k=\overline{T}^{t}_{\gamma}(\z^k)$, since
$\|\y^k\|\to+\infty$ and $\z^k\to\overline{\z}$. Therefore,
$\{\overline{T}^{t}_{\gamma}(\z^k)\}_{k\in\N}$ is bounded.
	

(iii) By passing to a subsequence, we may assume that $\cI(\z^k)\equiv \cI_0\subseteq[n]$, $\cB(\z^k)\equiv \cB_0$.
	Since $\overline{T}^t_{\gamma}(\z^k)_{\cB_0}=\z^k_{\cB_0}$, taking limits gives $\overline{\v}_{\cB_0} = \overline{\z}_{\cB_0}$. Moreover, since $\z^k \to \overline{\z}$, we have $\cB_0 \subseteq \cB(\overline{\z})$.
    It remains to handle the coordinates that are interior for every $z^k$ but
become boundary coordinates at the limit, i.e.,
$
    i\in\cB(\overline \z)\setminus \cB_0\subseteq \cI_0.
$
We show that $\overline{v}_i=\overline {z}_i$ for all such $i$.
	
	We proceed by contradiction. Suppose that there exists $i_0 \in \cB(\overline \z)\setminus \cB_0$ such that $\overline{v}_{i_0}\neq \overline{z}_{i_0}:=a $, and WLOG assume $a$ is the left endpoint of $\cl(\dom(\varphi))$. 
	Since $z^k_{i_0} \to a$ and $\overline{T}_\gamma^t(z^k)_{i_0} \to \overline{v}_{i_0} > a$, we have $z^k_{i_0} - \overline{T}_\gamma^t(\z^k)_{i_0} \to a - \overline{v}_{i_0} < 0$. Moreover, as $a$ is the left endpoint of $\cl(\dom(\varphi))$, it follows from Assumption \ref{assum:h} (iv) and  Definition~\ref{def:kernel} (ii) that $\varphi'(z^k_{i_0}) \to \varphi'(a) = -\infty$.  Hence, we have 
	\[
	\left(z^k_{i_0} - \overline{T}_\gamma^t(z^k)_{i_0}\right)\left(\varphi'(\overline{T}_\gamma^t(z^k)_{i_0}) - \varphi'(z^k_{i_0})\right) \to -\infty.
	\]
	Subsequently, by the convexity of $\varphi$, we have $(z^k_{i}-\overline{T}^{t}_{\gamma}(\z^{k})_{i})(\varphi^{\prime}(\overline{T}^{t}_{\gamma}(\z^{k})_i)-\varphi^{\prime}(z^{k}_i))\leq0$ for all $i\in[n]$, which leads to 
	\begin{equation}
		\label{eq:contra_blowup}
		\sum_{i\in\cI_0}\left(z^k_{i}-\overline{T}^{t}_{\gamma}(\z^{k})_{i}\right)\left(\varphi^{\prime}(\overline{T}^{t}_{\gamma}(\z^{k})_i)-\varphi^{\prime}(z^{k}_i)\right) \rightarrow -\infty. 
	\end{equation}
	On the other hand, we establish a finite lower bound for the same quantity in \eqref{eq:contra_blowup} by leveraging the optimality of $\overline{T}^t_\gamma(\z^k)$, which will lead to a contradiction. To this end, we define a sequence of interpolation points 
	\begin{equation}
		\label{eq:interpolate}
		{\z}^{\theta,k}:=  \theta\z^{k}+(1-\theta)\overline{T}^{t}_{\gamma}(\z^k),
	\end{equation}
	and a sequence of univariate functions as 
	\[{\phi}_k(\theta):=   G^{t}_{\gamma}\left({\z}^{\theta,k};\z^k\right),\qquad \forall\ \theta\in[0,1].
	\]
	By the optimality of $\overline{T}_\gamma^t(\z^k)$, we have $\phi_k(0) \leq \phi_k(\theta)$ for all $\theta \in [0,1]$, which implies
	$
	\frac{\phi_k(\theta) - \phi_k(0)}{\theta} \geq 0,
	$ for all $\theta \in (0,1]$.
	We next expand the difference quotient as follows: 
	\begin{equation*}
		\begin{aligned}
			\frac{\phi_k(\theta) -\phi_k(0)}{\theta} &=   \frac{\gamma\left(\z^{\theta,k};\z^k\right)-\gamma\left(\overline{T}^{t}_{\gamma}(\z^{k});\z^k\right)}{\theta}+\frac{g\left(\z^{\theta,k}\right)-g\left(\overline{T}^{t}_{\gamma}(\z^{k})\right)}{\theta} \\
			& \quad +\frac1t\sum_{i\in\cI_0}\frac{D_{\varphi}\left(z^{\theta,k}_i,z^k_i\right)-D_{\varphi}\left(\overline{T}^{t}_{\gamma}(\z^{k})_i,z^k_i\right)}{\theta}. 
		\end{aligned}
	\end{equation*}
	Letting $\theta \to 0_+$, and noting that $\z^{\theta,k} \to \overline{T}_\gamma^t(\z^k)$, we obtain for all $i\in\cI_0$, the following holds:
	\begin{align*}
		&\lim_{\theta \to 0_+} 
		\frac{D_{\varphi}\left(z^{\theta,k}_i, z^k_i\right) - D_{\varphi}\left(\overline{T}^{t}_{\gamma}(\z^{k})_i, z^k_i\right)}{\theta} \\
		= & \lim_{\theta \to 0_+} \left\{
		\frac{ \varphi\left(z^{\theta,k}_i\right) - \varphi\left(\overline{T}^{t}_{\gamma}(\z^{k})_i\right) }{\theta}
		- \left\langle \varphi'\left(z^k_i\right), \frac{z^{\theta,k}_i - \overline{T}^{t}_{\gamma}(\z^{k})_i}{\theta} \right\rangle \right\} \\
		= &\lim_{\theta \to 0_+} \left\{ \frac{ \varphi\left(\overline{T}^{t}_{\gamma}(\z^k)_i + \theta\left(z_i^k-\overline{T}^{t}_{\gamma}(\z^k)_i\right)\right) - \varphi\left(\overline{T}^{t}_{\gamma}(\z^{k})_i\right) }{\theta}\right\} - \varphi'(z^k_i) \left(z^{k}_i - \overline{T}^{t}_{\gamma}(\z^{k})_i\right)
		\\
		= & \left(z^k_i - \overline{T}^{t}_{\gamma}(\z^{k})_i\right) 
		\left( \varphi'\left(\overline{T}^{t}_{\gamma}(\z^{k})_i\right) - \varphi'(z^k_i) \right),
	\end{align*}
	where the first equality follows from the definition of Bregman divergence, the second one follows from \eqref{eq:interpolate}, and the last from the fact that $\varphi$ is continuously differentiable on $\intm(\dom(\varphi))$.  At boundary points, we use the extended directional derivative, i.e., $\varphi'(a) = -\infty$ and $\varphi'(c) = +\infty$.
	
	We now substitute the above limit into the expansion of the difference quotient, which yields
	\begin{align}
		\label{eq:diff_quo}
		\limsup_{\theta \to 0_+} \frac{\phi_k(\theta) - \phi_k(0)}{\theta}
		&= \limsup_{\theta \to 0_+} \left[
		\frac{ \gamma(\z^{\theta,k}; \z^k) - \gamma(\overline{T}_\gamma^t(\z^k); \z^k) }{\theta}
		+ \frac{ g(\z^{\theta,k}) - g(\overline{T}_\gamma^t(\z^k)) }{\theta} \right] \notag \\
		&\quad + \frac{1}{t} \sum_{i \in \mathcal{I}_0} (z^k_i - \overline{T}_\gamma^t(\z^k)_i)\left( \varphi'(\overline{T}_\gamma^t(\z^k)_i) - \varphi'(z^k_i) \right) \ge 0. 
	\end{align}
	Since the first two terms in the right-hand side of \eqref{eq:diff_quo} are uniformly bounded, as ensured by  the continuous differentiability of $\gamma(\cdot; \z^k)$, the local Lipschitz continuity of $g$, and the boundedness of the sequences $\{\z^k\}_{k \in \mathbb{N}}$ and $\{ \overline{T}^t_{\gamma}(\z^k) \}_{k \in \mathbb{N}}$, we obtain
	\[\liminf_{k\to\infty} \sum_{i\in\cI_0}\left(z^k_{i}-\overline{T}^{t}_{\gamma}(\z^{k})_{i}\right)\left(\varphi^{\prime}(\overline{T}^{t}_{\gamma}(\z^{k})_i)-\varphi^{\prime}(z^{k}_i)\right)>-\infty,\]
	which contradicts \eqref{eq:contra_blowup}. Therefore, we conclude that $\overline{\v}_{\cB(\overline{\z})}=\overline{\z}_{\cB(\overline{\z})}$, which directly implies $\cB(\overline{\z}) \subseteq \cB(\overline{\v})$.

Applying the same contradiction argument, now to a coordinate that belongs to
$\cB(\overline{\v})\setminus \cB(\overline{\z})$, we similarly obtain
$\overline{\v}_{\cB(\overline{\v})}
=
\overline{\z}_{\cB(\overline{\v})}$, and hence
$\cB(\overline{\v})\subseteq \cB(\overline{\z})$.
It follows that $\cB(\overline{\v})=\cB(\overline{\z})$, completing the proof.
\end{proof}

\section{Tightness of the Uniform Characterization in Theorem~\ref{th:Q}}
\label{sec:tight_th:Q}
    
The equivalence in Theorem~\ref{th:Q} is inherently uniform: The residual
must vanish along every interior feasible sequence converging to the boundary
point. In this appendix, we show that this uniform quantifier cannot, in
general, be weakened to the existence of a single such sequence. Specifically,
we construct an example in which there exists
$\{\z^k\}_{k\in\N}\subseteq \cX\cap\intm(\dom h)$ with
$\z^k\to \z$ and $R_\gamma^t(\z^k)\to0$, while the manifold stationarity
condition
$
    \bz\in \nabla_{\cM(\z)}f(\z)+\partial_{\cM(\z)}g(\z)
$
fails. 


\begin{example} \label{example:circle}
	 Consider $h(\x)=\sum_{i=1}^31/x_i^3$, $f(\x)=-10x_3$, and $g=\delta_{\cX}$ with 
	 \[\cX:=\textrm{conv}\left(\{(0,0,2)\}\cup(\cD\times\{1\})\right),\quad\text{where }\quad\cD=\left\{(x_1,x_2)\in\R^2_+:x_1^2+(x_2-1)^2\leq1\right\}.\]
	\end{example}
Let $t=1$, $\gamma(\y;\x)=-10y_3$, and $\overline{\z}:=(0,0,1)$. Then, we have 
\[
    \cM(\overline{\z})=\{\x\in\mathbb R^3:x_1=x_2=0\},
    \quad \text{and} \quad 
    \cM(\overline{\z})\cap\cX=\{(0,0,x_3):1\le x_3\le2\}.
\]
Hence, we know $
    T_{\cM(\overline{\z})}=\{(0,0,s):s\in\mathbb R\}.
$
Since $g=\delta_{\cX}$, the Riemannian subdifferential of $g$ at $\overline{\z}$
relative to $\cM(\overline{\z})$ is
\[
    \partial_{\cM(\overline{\z})}g(\overline{\z})
    =
    \{(0,0,s):s\le0\}.
\]
Moreover,
$
    \nabla_{\cM(\overline{\z})}f(\overline{\z})=(0,0,-10).
$
Therefore,
\[
    \bz\notin \nabla_{\cM(\overline{\z})}f(\overline{\z})+\partial_{\cM(\overline{\z})}g(\overline{\z}).
\]


Next, we show that the residual can vanish along a particular interior feasible
path approaching \(\overline{\z}\). Let
\[
    \z^k:=\left(\sin(\theta_k),\,1-\cos(\theta_k),\,1\right),
    \qquad
    \theta_k\in\left(0,\frac{\pi}{3}\right],\quad \theta_k\downarrow0 .
\]
Since
$    \sin^2(\theta_k)+\bigl((1-\cos(\theta_k))-1\bigr)^2=1
$
and both $\sin(\theta_k)$ and $1-\cos(\theta_k)$ are positive, we have
$\z^k\in\cX\cap\intm(\dom h)$. Moreover, $\z^k\to\overline{\z}$.
We will prove that
$
    R_\gamma^t(\z^k)\to0 .
$
We first collect a few elementary estimates needed for the argument. 

\begin{lemma}\label{le:circle}
    Let $\cD_{\alpha}:=\{(x_1,x_2):\x\in\cX,~x_3\geq1+\alpha\}$ for $\alpha\in(0,1)$. The following hold:
	\begin{enumerate}[{\rm (i)}]
		\item $\cD_{\alpha}=(1-\alpha)\cdot \cD$.
		\item $\dist((z^k_1,z^k_2),\cD_{\alpha})\geq \alpha z^k_2$ and $ z^k_2\geq \frac{1}{2}(z^k_1)^2$ for $k\in\N$.
		\item $\lim_{x\to0+}D_{\varphi}(x+bx^2,x)=+\infty$ for any $b\in\R\setminus\{0\}$.
	\end{enumerate}
\end{lemma}
Suppose that $T^t_{\gamma}(\z^k)_3\geq 1+\alpha$ for some $\alpha\in(0,1)$ for all $k\in\N$. Then, we have $(T^t_{\gamma}(\z^k)_1,T^t_{\gamma}(\z^k)_2)\in\cD_{\alpha}$ by the definition of $\cD_{\alpha}$, and
Lemma \ref{le:circle} (ii) implies  
\[\left\|(z^k_1,z^k_2)-\left(T^t_{\gamma}(\z^k)_1,T^t_{\gamma}(\z^k)_2\right) \right\|\geq \alpha z^k_2.  \]
Then, there exists $i_k\in\{1,2\}$ such that
\[\left|z^k_{i_k}-T^t_{\gamma}(\z^k)_{i_k}\right|\geq\frac{\sqrt{2}\alpha}{2} z^k_2.  \]
Observe that $z^k_2\geq(z^k_2)^2$ by $z^k_2\leq1$, and $z^k_2\geq \frac12(z^k_1)^2$ by Lemma \ref{le:circle} (ii). We further have
\[\left|z^k_{i_k}-T^t_{\gamma}(\z^k)_{i_k}\right|\geq\frac{\sqrt{2}\alpha}{4} (z^k_{i_k})^2.  \]
By the one-dimensional monotonicity property of Bregman divergences, for fixed $y$, namely
\[
    D_\varphi(z,y)\ge D_\varphi(x,y)
    \quad
    \text{whenever either } z\le x\le y
    \text{ or } y\leq x\leq z,
\]
we see that \[D_{\varphi}\left(T^t_{\gamma}(\z^k)_{i_k},z^k_{i_k}\right)\geq \min\left(D_{\varphi}\left( z^k_{i_k}-\frac{\sqrt{2}\alpha}{4} (z^k_{i_k})^2,z^k_{i_k}\right),D_{\varphi}\left( z^k_{i_k}+\frac{\sqrt{2}\alpha}{4} (z^k_{i_k})^2,z^k_{i_k}\right)\right).\]
Together with Lemma~\ref{le:circle}(iii) and \(z^k_{i_k}\to0+\), this implies
\[
    D_h\left(T^t_{\gamma}(\z^k),\z^k\right)
    \ge
    D_{\varphi}\left(T^t_{\gamma}(\z^k)_{i_k},z^k_{i_k}\right)
    \to+\infty .
\]
On the other hand, by the optimality of $T^t_{\gamma}(\z^k)$, we have
\begin{equation}\label{eq:circleG}
    -10T^t_{\gamma}(\z^k)_3
    +
    D_h\left(T^t_{\gamma}(\z^k),\z^k\right)
    \le
    G^t_{\gamma}\left(\z^k;\z^k\right)
    =
    -10z^k_3 .
\end{equation}
Since $z_3^k=1$ and $T^t_{\gamma}(\z^k)\in\cX$, we have
$1\le T^t_{\gamma}(\z^k)_3\le2$, and hence \eqref{eq:circleG} gives
\[
    D_h\left(T^t_{\gamma}(\z^k),\z^k\right)
    \le
    10\left(T^t_{\gamma}(\z^k)_3-1\right)
    \le 10,
\]
a contradiction. Therefore, for every $\alpha\in(0,1)$, the inequality
$T^t_{\gamma}(\z^k)_3\ge1+\alpha$ cannot hold along any subsequence. Since
$T^t_{\gamma}(\z^k)_3\ge1$, we conclude that
$    T^t_{\gamma}(\z^k)_3\to1.
$ Using \eqref{eq:circleG} once more, we obtain
\[
    0\le
    D_h\left(T^t_{\gamma}(\z^k),\z^k\right)
    \le
    10\left(T^t_{\gamma}(\z^k)_3-1\right)\to0.
\]
Thus
\[
    R^t_{\gamma}(\z^k)
    =
    D_h\left(T^t_{\gamma}(\z^k),\z^k\right)
    \to0,
\]
as claimed.



It remains to prove Lemma~\ref{le:circle}. 
\begin{proof}[Proof of Lemma \ref{le:circle}]
	(i) By the definition of $\cX$, we have
	\[ \cX=\left\{\left((1-\tau)x_1,(1-\tau)x_2,1+\tau\right):\tau\in[0,1],~(x_1,x_2)\in\cD\right\}.\]
Together with the definition of $\cD_\alpha$, this gives
\begin{equation}\label{eq:Dalpha}
\begin{aligned}
	\cD_{\alpha}&=\left\{\left((1-\tau)x_1,(1-\tau)x_2\right):\alpha\leq\tau\leq1,~(x_1,x_2)\in\cD\right\}\\
	&=\bigcup_{\tau\in[\alpha,1]}(1-\tau)\cdot\cD\supseteq(1-\alpha)\cdot \cD. \\ 
\end{aligned}   
\end{equation}
On the other hand, by $\bz\in\cD$ and the convexity of $\cD$, for all $\tau\in[\alpha,1]$,
\[(1-\tau)\cdot\cD=\frac{1-\tau}{1-\alpha}\cdot(1-\alpha)\cD+\left(1-\frac{1-\tau}{1-\alpha} \right)\cdot\bz\subseteq(1-\alpha)\cD.  \]
Combined with \eqref{eq:Dalpha}, it yields $\cD_{\alpha}=(1-\alpha)\cD$ as desired.

(ii) 
By part (i), we have $\cD_\alpha=(1-\alpha)\cD$. Hence
\[
    \cD_\alpha
    \subseteq
    \left\{
        (u_1,u_2)\in\R^2:
        u_1^2+(u_2-(1-\alpha))^2\le (1-\alpha)^2
    \right\}.
\]
Therefore,
\[
\begin{aligned}
\dist((z^k_1,z^k_2),\cD_\alpha)
&\ge
\sqrt{\sin^2(\theta_k)+(\alpha-\cos(\theta_k))^2}-(1-\alpha)  \\
&=
\sqrt{1+\alpha^2-2\alpha\cos(\theta_k)}-(1-\alpha)\\
&=
\frac{2\alpha(1-\cos(\theta_k))}
{\sqrt{1+\alpha^2-2\alpha\cos(\theta_k)}+(1-\alpha)}\\
&\ge
\alpha(1-\cos(\theta_k))
=
\alpha z^k_2 .
\end{aligned}
\]
Here, the last inequality follows from \(\theta_k\in(0,\pi/3]\), which gives
\(\cos(\theta_k)\ge 1/2\), and hence
\[
    1+\alpha^2-2\alpha\cos(\theta_k)
    \le 1+\alpha^2-\alpha
    \le 1.
\]
Thus the denominator above is at most \(2\).
Moreover,
\[
z^k_2
=
1-\cos(\theta_k)
=
2\sin^2\left(\frac{\theta_k}{2}\right)
\ge
2\sin^2\left(\frac{\theta_k}{2}\right)
\cos^2\left(\frac{\theta_k}{2}\right)
=
\frac12\sin^2(\theta_k)
=
\frac12(z^k_1)^2 .
\]


(iii) For sufficiently small $x>0$, we have $x+bx^2>0$.  Direct computation yields
\[ D_{\varphi}(x+bx^2,x)=\frac{1}{(x+bx^2)^3}-\frac{1}{x^3}+\frac{3bx^2}{x^4}=\frac{1}{x^3}\left(\frac{1}{(1+bx)^3}-1+3bx\right)=\frac{6b^2+o(1)}{x(1+bx)^3}. \]
Since $b\neq0$ and $1+bx\to1$, the last expression diverges to $+\infty$
as $x\to0+$.
\end{proof}
\end{document}